\documentclass[11p, reqno]{amsart}
\usepackage{amsfonts}
\usepackage{amsmath}
\usepackage{amssymb}
\usepackage[pdftex]{graphicx,color}

\newtheorem{theorem}{Theorem}[section]
\newtheorem{lemma}[theorem]{Lemma}

\newtheorem{proposition}[theorem]{Proposition}
\newtheorem{remark}[theorem]{Remark}

\numberwithin{equation}{section}

\DeclareMathOperator{\re}{Re}
\DeclareMathOperator{\im}{Im}
\usepackage{graphicx}
\graphicspath{%
    {converted_graphics/}
    {/}
}
\begin{document}
\def\R{{\mathbb R}}
\def\C{{\mathbb C}}
\def\N{{\mathbb N}}
\def\S{{\mathbb S}}
\def\DD{{\mathbb D}}
\def\rr{{\cal R}}
\def\e{\emptyset}
\def\dQ{\partial Q}
\def\dk{\partial K}
\def\endofproof{{\rule{6pt}{6pt}}}
\def\di{\displaystyle}
\def\dist{\mbox{\rm dist}}
\def\u-{\overline{u}}
\def\du{\frac{\partial}{\partial u}}
\def\dv{\frac{\partial}{\partial v}}
\def\dt{\frac{d}{d t}}
\def\dx{\frac{\partial}{\partial x}}
\def\con{\mbox{\rm const }}
\def\Box{\spadesuit}
\def\ii{{\bf i}}
\def\curl{{\rm curl}\,}
\def\dive{{\rm div}\,}
\def\grad{{\rm grad}\,}
\def\dist{\mbox{\rm dist}}
\def\pr{\mbox{\rm pr}}
\def\pp{{\cal P}}
\def\supp{\mbox{\rm supp}}
\def\Arg{\mbox{\rm Arg}}
\def\In{\mbox{\rm Int}}
\def\Re{\mbox{\rm Re}\:}
\def\li{\mbox{\rm li}} 
\def\ep{\epsilon}
\def\tr{\tilde{R}}
\def\be{\begin{equation}}
\def\ee{\end{equation}}
\def\cn{{\mathcal N}}
\def\sn{{\mathbb  S}^{n-1}}
\def\Ker {{\rm Ker}\:}
\def\el{E_{\lambda}}
\def\Rc{{\mathcal R}}
\def\Ha{H_0^{ac}}
\def\la{\langle}
\def\ra{\rangle}
\def\Ko{\Ker G_0}
\def\Kd{\Ker G_b}
\def\hc{{\mathcal H}}
\def\caH{\hc}
\def\caO{{\mathcal O}}
\def\la{\langle}
\def\ra{\rangle}
\def\sp{\sigma_{+}}
\def\et{|\xi'|^2}
\def\pa{\partial}
\def\ep{\epsilon}
\def\cab{C_{\alpha, \beta}}
\def\lg{L^2(\Gamma)}
\def\tu{\tilde{u}}
\def\ep{\epsilon}
\def\oc{{\mathcal O}}
\def\ac{{\mathcal A}}
\def\nc{{\mathcal N}}

\title[Location of eigenvalues] {Location of eigenvalues for the wave equation with dissipative boundary conditions}
\author[V. Petkov]{Vesselin Petkov}
\address {Universit\'e de Bordeaux, Institut de Math\'ematiques de Bordeaux,  351, Cours de la Lib\'eration, 33405  Talence, France}

\email{petkov@math.u-bordeaux.fr}


\thanks{The author was partially supported by the ANR project Nosevol BS01019 01}


\subjclass[2000]{Primary 35P20, Secondary 47A40, 35L05}

\date{}

\dedicatory{}


\begin{abstract} We examine the location of the eigenvalues of the generator $G$ of a semi-group $V(t) = e^{tG},\: t \geq 0,$ related to the wave equation in an unbounded domain $\Omega \subset \R^d$ with dissipative boundary condition $\partial_{\nu}u - \gamma(x) \pa_t u = 0$ on $\Gamma = \partial \Omega.$ We study two cases: $(A): \: 0 < \gamma(x) < 1,\: \forall x \in \Gamma$ and $(B):\: 1 < \gamma(x), \: \forall x \in \Gamma.$ We prove that for every $0 < \ep \ll 1,$ the eigenvalues of $G$ in the case $(A)$ lie in the region
$\Lambda_{\epsilon} = \{ z \in \C:\: |\re z | \leq C_{\epsilon} (|\im z|^{\frac{1}{2} + \epsilon} + 1), \: \re z < 0\},$
while in the case $(B)$ for every $0 < \ep \ll 1$ and every $N \in \N$ the eigenvalues lie in $\Lambda_{\ep} \cup {\mathcal R}_N,$ where
${\mathcal R}_N = \{z \in \C:\: |\im z| \leq C_N (|\re z| + 1)^{-N},\: \re z < 0\}.$   
\end{abstract}

\maketitle

\section{Introduction}
 
Let $K \subset \R^d,$ $d \geq 2$, be an open bounded domain and let $\Omega = \R^d \setminus \bar{K}$ be connected. We suppose that the boundary $\Gamma$ of $\Omega$ is $C^{\infty}.$ 
Consider the boundary problem
\begin{equation} \label{eq:1.1}
\begin{cases} u_{tt} - \Delta u = 0 \: {\rm in}\: \R_t^+ \times \Omega,\\
\partial_{\nu}u - \gamma(x) \pa_t u= 0 \: {\rm on} \: \R_t^+ \times \Gamma,\\
u(0, x) = f_1, \: u_t(0, x) = f_2 \end{cases}
\end{equation}
with initial data $(f_1, f_2) \in H^1(\Omega) \times L^2(\Omega) = {\mathcal H}.$
Here $\nu(x)$ is the unit outward normal at $x \in \Gamma$ pointing into $\Omega$ and $\gamma(x)\geq 0$ is a $C^{\infty}$ function on $\Gamma.$ The solution of the problem (\ref{eq:1.1}) is given by $(u(t, x) , u_t(t, x))= V(t)f = e^{tG} f,\: t \geq 0$, where $V(t)$ is a contraction semi-group in ${\mathcal H}$ whose  generator
$$ G = \Bigl(\begin{matrix} 0 & 1\\ \Delta & 0 \end{matrix} \Bigr)$$
has a domain $D(G)$ which is the closure in the graph norm of functions $(f_1, f_2) \in C_{(0)}^{\infty}(\bar{\Omega}) \times C_{(0)}^{\infty} (\bar{\Omega})$ satisfying the boundary condition $\partial_{\nu} f_1 - \gamma f_2 = 0$ on $\Gamma.$ For $d$ odd Lax and Phillips \cite{LP} proved that the spectrum of $G$ in $\Re z < 0$ is formed by isolated eigenvalues with finite multiplicity, while the continuous spectrum of $G$ coincides with $\ii R.$  We obtain the same result for all dimensions $d \geq 2$ under the restriction $\gamma(x) \neq 1$ in the case $d$ even by using the Dirichlet-to-Neumann map $\nc(\lambda)$ (see Section 6). 
Notice that if $Gf =\lambda f$ with $f = (f_1, f_2) \neq 0, \: \re \lambda < 0$ and $\partial_{\nu} f_1 - \gamma f_2 =0$ on $\Gamma$, we get 
\begin{equation} \label{eq:1.2}
\begin{cases} (\Delta - \lambda^2) f_1 = 0 \:{\rm in}\: \Omega,\\
\partial_{\nu} f_1 -  \lambda \gamma f_1 = 0\: {\rm on}\: \Gamma \end{cases}
\end{equation}
and $V(t) f =  e^{\lambda t}f  $ has an exponentially decreasing global energy. Such solutions are called {\bf asymptotically disappearing} and they perturb the inverse scattering problems. Recently it was proved \cite{CPR} that if we have at least one eigenvalue $
\lambda$ of $G$ with $\re \lambda < 0$, then the wave operators $W_{\pm}$ related to the Cauchy problem for the wave equation and the boundary problem (\ref{eq:1.1}) are not complete, that is ${\text Ran}\: W_{-} \not=  {\text Ran}\: W_{+}$. Hence  we cannot define the scattering operator $S$  by the product $ W_{+}^{-1}\circ W_{-}$. Notice that if the global energy is conserved in time and the unperturbed and perturbed  problems are associated to unitary groups, the corresponding scattering operator $S(z): L^2(\S^{d-1}) \to L^2(\S^{d-1})$ satisfies the identity
$$S^{-1}(z)= S^*(\bar{z}),\: z \in \C,$$
if $S(z)$ is invertible at $z$.  Since $S(z)$ and $S^*(z)$ are analytic operator-valued operators in the "physical" half plane $\{z \in \C:\im z < 0\}$ (see \cite{LP1}) the above relation implies that $S(z)$ is invertible for $\im z > 0$. For dissipative boundary problems this relation in general is not true and $S(z_0)$ may have a non trivial kernel for some $z_0, \im z_0 > 0.$ For odd dimensions $d$ Lax and Phillips \cite{LP} proved that this implies that $\ii z_0$ is an eigenvalue of $G$. Thus the analysis of the location of the eigenvalues of $G$ is important for the location of the points where  the kernel of $S(z)$ is not trivial.\\

In the scattering theory of Lax-Phillips \cite{LP} for odd dimensions $d$ the energy space can be presented as a direct sum
${\mathcal H} = D_a^{-} \oplus K_a \oplus D_a^{+},\: a > 0,$
and we have the relations
$$V(t) D_a^{+} \subset D_a^{+},\: V(t) (K_a) \subset K_a \oplus D_{+}^a,\: V(t) D_a^{-} \subset {\mathcal H},\: t \geq 0.$$
R. Phillips defined a system as non controllable if there exists a state $f \in K_a$ such that
$V(t) f \perp D_a^{+},\: t \geq 0.$
This means that there exists states in the "black box" $K_a$ which remain undetected  by the scattering process. Majda \cite{M} proved that if we have such state $f$, then  $(u(t, x), u_t(t, x)) = V(t) f$ is a disappearing solution, that is there exists $T > 0$ depending on $f$ such that $u(t, x)$  vanishes for all $t \geq T > 0$. On the other hand, if $\gamma(x) \neq 1,  \forall x \in \Gamma,$ and the boundary is analytic there are no disappearing solutions (see \cite{M}). Thus in this case  it is natural to search asymptotically disappearing solutions. The existence of examples in the case
$\gamma \equiv 1$ when the point spectrum of $G$ is empty has been mentioned in \cite{M1}. Since we did not found  a proof of this result in the literature, for reader convenience we propose a simple analysis  of this question for the ball $B_3 = \{ x \in \R^3:\: |x| \leq 1\}$. In the  Appendix we prove that if $\gamma \equiv 1$ and $\bar{K} = B_3$, the generator $G$ has no eigenvalues in $\{z \in \C:\Re z < 0\}.$ \\

We study in the Appendix also the case when $\gamma \equiv {\rm  const}\: \neq 1$ and $K = B_3.$ If $0 < \gamma < 1$, we show that there are no real eigenvalues of $G$ and the complex eigenvalues 
$\lambda$ lie in the region 
$$\{\lambda \in \C: \: |\re \lambda| < |\im \lambda|, \: \re \lambda < 0\}.$$
 On the other hand, for $\gamma > 1$ all eigenvalues of $G$ are real and they lie in the interval $ [-\frac{1}{\gamma-1}, - \infty)$. Moreover, in this case there are infinite number real  eigenvalues of $G$ and when $\gamma \searrow 1$ the eigenvalues of $G$ go to $-\infty$. For arbitrary strictly convex obstacle $K$ and $\gamma(x) > 1,\: \forall x\in \Gamma,$ we obtain a similar result in Theorem 1.3 proving that with exception of a finite number eigenvalues all other are confined in a  very small neighbourhood of the negative real axis.\\

If $\max_{x \in \Gamma}|\gamma(x)- 1|$ is sufficiently small, the leading term of the back-scattering amplitude $a(\lambda, -\omega, \omega),\: \omega \in \sn,$ becomes very small for all directions $\omega \in \sn$ and for $\gamma \equiv 1$ this leading term vanishes for all directions (see \cite{M2}).  For strictly convex obstacles and $\gamma \equiv 1$ the second term of the back-scattering amplitude does not vanish (see \cite{GA}), but it is negligible for the applications  (see Section 5 in \cite{P} for the case of first order systems).  The existence of a space with infinite dimension of eigenfunctions of $G$ implies that one has a large set of initial data for which the solutions of (\ref{eq:1.1}) are asymptotically disappearing. Notice that these solutions cannot be outgoing in the sense of Lax-Phillips (see \cite{M}), that is they have a non-vanishing projection on the space $D_a^{-}$ mentioned above. Moreover, the  eigenvalues of $G$ are stable under perturbations of the boundary and the boundary condition (see \cite{CPR}).\\ 

Now we pass to the description of our results.
In \cite{M1} Majda examined the location of the eigenvalues of $G$ and he proved that
 if $\sup \gamma(x) < 1,$  the eigenvalues of $G$ lie in the region
$$E_1 = \{ z \in \C: \: |\re z | \leq C_1 (|\im z|^{3/4} + 1),\: \re z < 0\},$$
while if $\sup \gamma(x) \geq 1,$ the eigenvalues of $G$ lie in $E_1 \cup E_2$, where
$$E_2 = \{ z \in \C: \: |\im z| \leq C_2 (|\re z|^{1/2} + 1),\: \re z < 0\}.$$

The purpose of this paper is to improve the above results for the location of eigenvalues. We consider two cases:
$(A): \: 0 < \gamma(x) < 1,\: \forall x \in \Gamma,\:(B): \: \gamma(x) > 1, \:\forall x \in \Gamma.$ Our main result is the following
\begin{theorem} In the case $(A)$ for every $\epsilon, \: 0 < \epsilon \ll 1,$ the eigenvalues of $G$ lie in the region
$$\Lambda_{\epsilon} = \{ z \in \C:\: |\re z | \leq C_{\epsilon} (|\im z|^{\frac{1}{2} + \epsilon} + 1), \: \re z < 0\}.$$
In the case $(B)$ for every  $\epsilon, \: 0 < \epsilon \ll 1,$ and every $N \in \N$ the eigenvalues of $G$ lie in the region $\Lambda_{\epsilon} \cup {\mathcal R}_N$, where
$${\mathcal R}_N = \{z \in \C:\: |\im z| \leq C_N (|\re z| + 1)^{-N},\: \re z < 0\}.$$
\end{theorem}

For strictly convex obstacles $K$ we prove a better result in the case (B). 
\begin{theorem} Assume that $K$ is strictly convex. 
In the case $(B)$ there exists $R_0 > 0$ such that for  every $N \in \N$ the eigenvalues of $G$ lie in the region $\{ z \in \C: |z| \leq R_0,\: \re z < 0\} \cup {\mathcal R}_{N}.$ 
\end{theorem}

The eigenvalues of $G$ are symmetric with respect to the real axis, so it is sufficient to examine the location of the eigenvalues whose imaginary part is  non negative.
Introduce in $\{ z \in \C: \: \im z \geq 0\}$ the sets
$$Z_1 = \{ z\in \C:\: \re z = 1,  h^{\delta} \leq \im z \leq 1\}, \:0 < h \ll 1, \: 0 < \delta < 1/2,$$
$$ Z_2 = \{ z \in \C: \re z =  - 1, 0 \leq \im z \leq 1\},\: Z_3 = \{ z \in \C: |\re z| \leq 1,\: \im z = 1\}.$$

\begin{figure}[tbp] 
  \centering
  \includegraphics[bb=0 0 667 402,width=5.57in,height=2.90in,keepaspectratio]{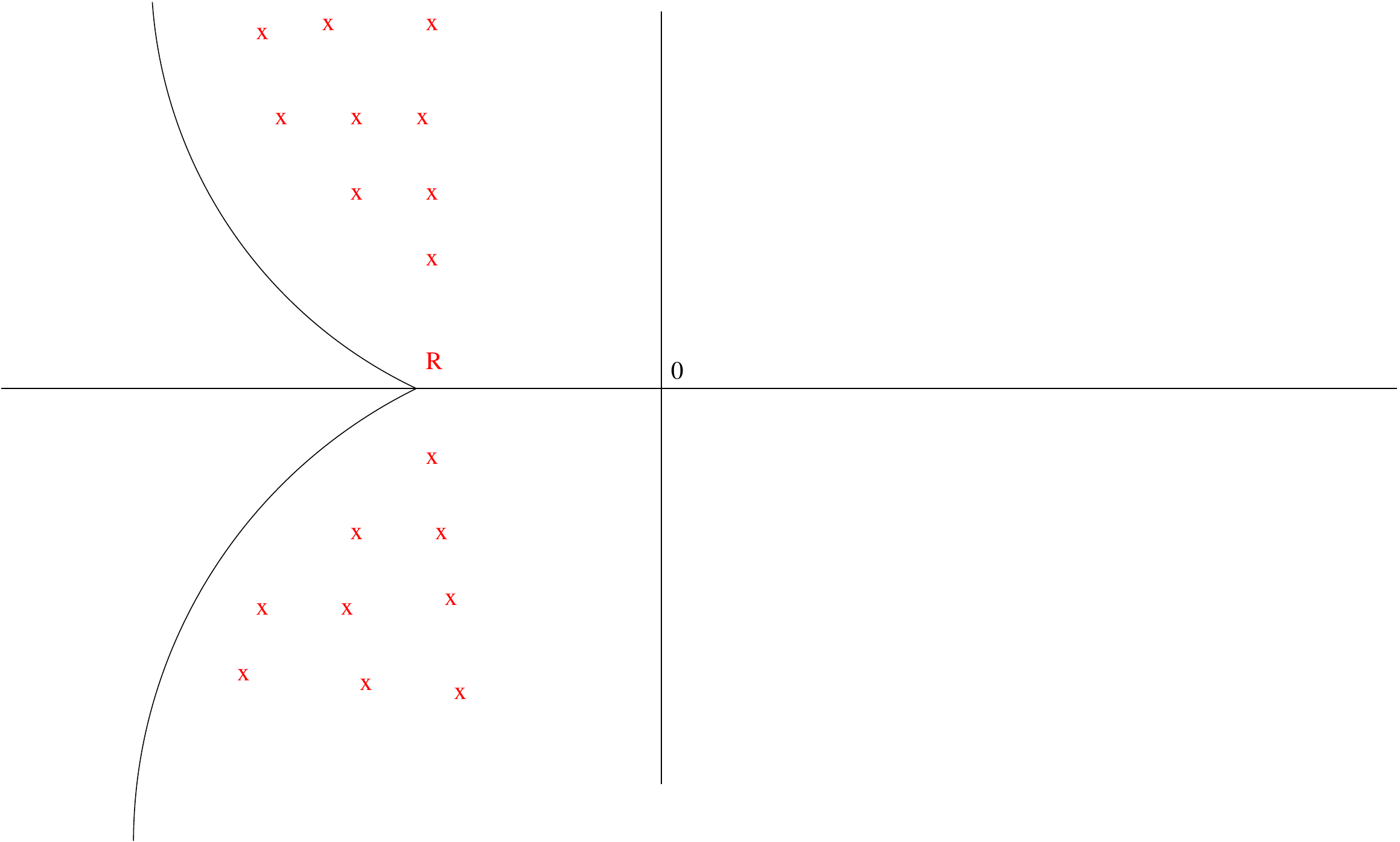}
  \caption{Eigenvalues, $0 < \gamma(x) < 1$}
  \label{fig:values3}
\end{figure}

\begin{figure}[tbp] 
  \centering
  \includegraphics[bb=0 0 487 402,width=4.67in,height=2.68in,keepaspectratio]{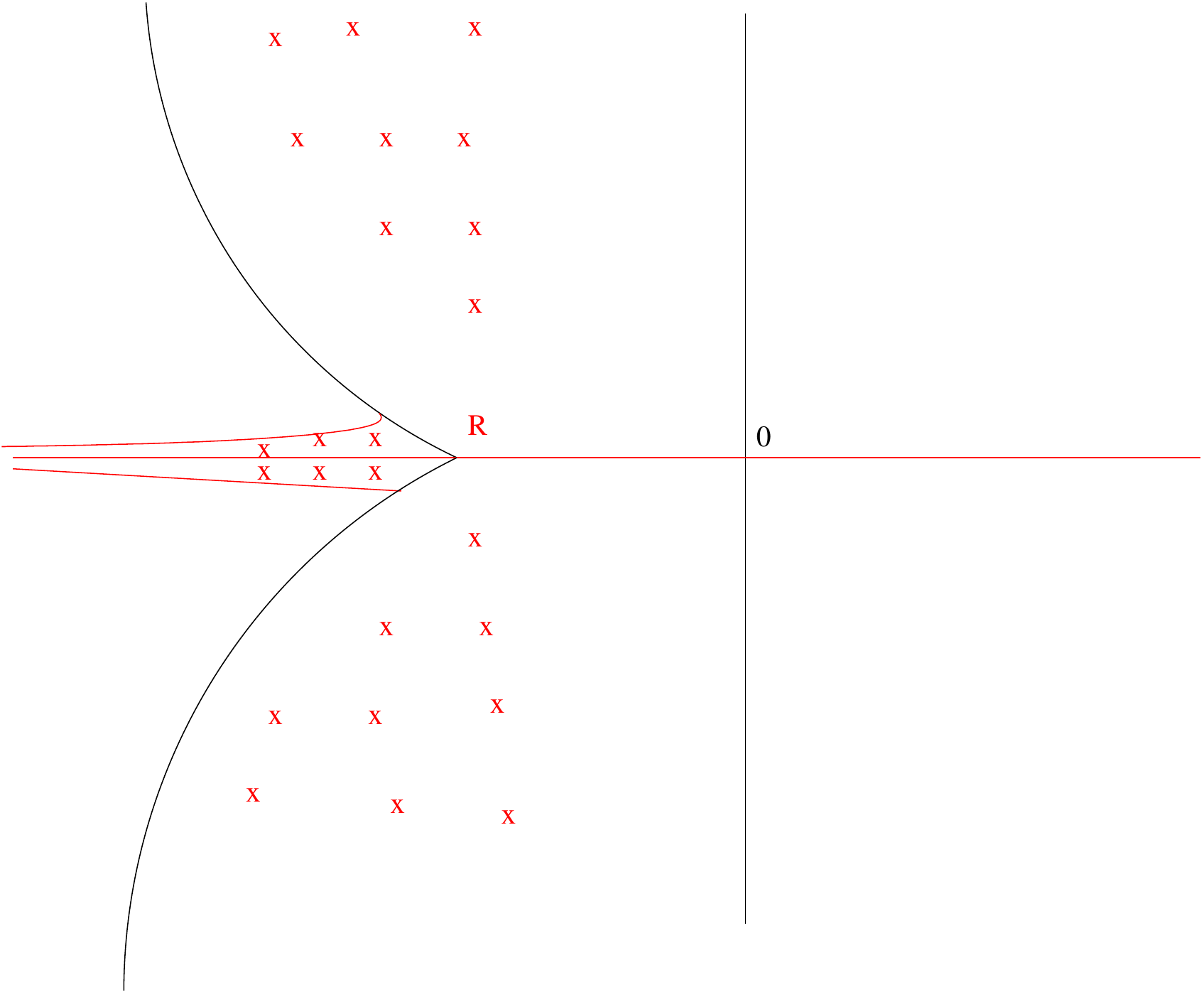}
  \caption{Eigenvalues, $1 < \gamma(x)$}
  \label{fig:values5}
\end{figure}

We put $\lambda = \frac{ \ii \sqrt{z}}{h}$ and we use the branch $0 \leq \arg z < 2 \pi$ with $\im \sqrt{z} > 0$ if $\im z > 0.$ From (\ref{eq:1.2}) we deduce that the eigenfunctions $u$ of $G$ satisfy the problem
\begin{equation} \label{eq:1.3}
\begin{cases} (-h^2\Delta - z) u = 0 \:{\rm in}\: \Omega,\\
- \ii h\partial_{\nu}u -   \gamma \sqrt{z}u = 0\: {\rm on}\: \Gamma. \end{cases}
\end{equation}


The proofs of Theorems 1.2 and 1.3 are based on a semi-classical analysis of the equation
\begin{equation} \label{eq:1.4}
 N_{ext}(z, h)f - \gamma \sqrt{z} f = 0,
\end{equation}
where $f = u\vert_{x \in \Gamma}$ is the trace of an eigenfunction of $G$. 
Here 
$$\nc_{ext}(z, h)f:\:H_h^s(\Gamma) \ni f \longrightarrow h D_{\nu} u\vert_{\Gamma} \in H_h^{s-1}(\Gamma)$$ 
is the exterior Dirichlet-to-Neumann map, 
 $D_{\nu} = -\ii \pa_{\nu}$ and $u$ is the solution of the problem
\begin{equation} \label{eq:1.5}
\begin{cases} (-h^2 \Delta - z) u = 0\: {\rm in}\: \Omega,\: u \in H_h^2(\Omega),\\
u\vert_{x \in \Gamma} = f.\end{cases}
\end{equation}
In the paper we use the semi-classical Sobolev space $H_h^s(\Gamma),\: s \in \R,$ with norm $\|\langle h D \rangle^s u\|_{\lg}$,
 where $\langle h D \rangle = (1 + (h D_x)^2)^{1/2}.$ The purpose is to prove that if $z \in Z_1 \cup Z_2 \cup Z_3$ lies in some regions and $h$ is small enough from (\ref{eq:1.4}) we get $f = 0$ which is not possible for an eigenfunction $u$. In this direction our strategy is close to that for the analysis of eigenvalues-free regions for the interior transmission eigenvalues in \cite{V} and \cite{V1}. We apply some results for the interior Dirichlet-to-Neumann map $\nc_{int}(z, h)$  established in \cite{V} and \cite{V1} for bounded domains which after modifications and some constructions remain true for the exterior Dirichlet-to-Neumann map $\nc_{ext}(z, h)$ defined above. \\
 
 The paper is organized as follows. In Section 2 we collect some results concerning the semi-classical exterior Dirichlet-to-Neumann map $\nc_{ext}(z, h).$ The eigenvalues-free regions for the case (A) are discussed in Section 3. In Section 4 we study the case (B), where the arguments for the case (A) are not applicable for the investigation of eigenvalues close to the negative real axis. The strictly convex obstacles are examined in Section 5. In Section 6 we discuss the question of the discreteness of the spectrum of $G$ in $\{z\in \C:\: \re z < 0\}$ for dimensions $d \geq 2.$ For odd dimensions $d$, as it was mentioned above, this result was obtained in the classical paper \cite{LP}. For $d$ even we present a proof based on the properties of the Dirichlet-to-Neumann map $\nc_{ext}(\lambda)$. Moreover, we obtain a trace formula for the counting function of the eigenvalues of $G$ in an open domain  $\omega \subset \{z\in \C:\: \re z < 0\}$. Finally, in the Appendix we examine the special case when $K$ is unit ball in $\R^3$ and $\gamma$ is a constant.

\section{Dirichlet-to-Neumann map}

 In our exposition we apply some $h$-pseudo-differential operators and we are going to recall some basic facts. Let $X$ be a $C^{\infty}$ smooth compact manifold without
boundary with dimension $d-1 \geq 1$. Let $(x, \xi)$ be the canonical local coordinates in $T^*(X)$ and let $a(x, \xi, h) \in C^{\infty}(T^*(X)).$ 
 Given $m \in \R, \: l \in \R, \delta > 0$ and a function $c(h) > 0$, one denotes by $S_{\delta}^{l, m}(c(h))$ the set of symbols $a(x, \xi, h)$ such that
$$ |\pa_{x}^{\alpha} \pa_{\xi}^{\beta} a(x, \xi, h)| \leq C_{\alpha, \beta} (c(h))^{l - \delta(|\alpha| + |\beta|)} (1 + |\xi|)^{m - |\beta|},\: (x, \xi) \in T^*(X),\:\forall \alpha, \forall \beta.$$
If $c(h) = 1$, we denote $S_{\delta}^{l, m}(c(h))$ simply by $S_{\delta}^{l, m}$ and the symbols restricted to a domain where $|\xi|\leq C$ will be denoted by $a \in S_{\delta}^l(c(h)).$  We use also symbols $a(x, \xi, h) \in S_{0, 1}^m$ satisfying the estimates
$$ |\pa_{x}^{\alpha} \pa_{\xi}^{\beta} a(x, \xi, h)| \leq C_{\alpha, \beta}  (1 + |\xi|)^{m - |\beta|},\: (x, \xi) \in T^*(X),\:\forall \alpha, \forall \beta.$$

One defines the $h-$pseudo-differential operator $Op_h(a)$ with symbol $a(x, \xi, h)$ by
$$(Op_h(a) f)(x) = (2 \pi h)^{-d+1}\int_{T^*X} e^{-\ii \langle x - y, \xi \rangle /h} a(x, \xi, h) f(y) dy d \xi.$$

For the reader convenience we recall two properties of the semi-classical pseudo- differential operators $Op_h(a)$ (see Section 7 of \cite{DS} and  Proposition 2.1 of \cite{V}). Assume that $a \in C^{\infty}(T^*(X))$ satisfies the bounds
\begin{equation} \label{eq:2.1}
|\pa_x^{\alpha} a(x, \xi, h)| \leq c_0(h) h^{-|\alpha|/2}, \: (x, \xi) \in T^*(X)
\end{equation}
for $|\alpha| \leq d+1$, where $c_0(h) > 0$ is a parameter. Then there exists a constant $C > 0$ independent of $h$ such that
\begin{equation} \label{eq:2.2}
\|Op_h(a)\|_{L^2(X) \to L^2(X)} \leq C c_0(h).
\end{equation}

Next for $0 \leq \delta < 1/2$ we have a calculus and if $a \in S_{\delta}^{l_1, m_1}, \: b \in S_{\delta}^{l_2, m_2}$, then for $s \in \R$ we get
$$\|Op_h(a) Op_h(b) - Op_h (a b) \|_{H^s(X) \to H^{s - m_1 - m_2+1}(X)} \leq C h^{-l_1 - l_2 - 2\delta +1}.$$
We refer to \cite{DS} for more details concerning the calculus.
The left hand side of last inequality can be estimated also in some cases when one  of the symbols $a$ or $b$ is in a class $S_{\delta}^{l, m}$ with
$0 \leq  \delta < 1$. For the  precise statements the reader should consult Proposition 2.2 in \cite{V} and Proposition 4.2 in \cite{V1}.\\

Let $(x', \xi')$ be the coordinates on $T^*(\Gamma)$. Denote by $r_0(x', \xi')$ the principal symbol of the Laplace-Beltrami operator $-\Delta_{\Gamma}$ on $\Gamma$ equipped with the Riemannian metric induced by the Euclidean metric in $\R^d.$
For $ z \in Z_1 \cup Z_2 \cup Z_3$ let 
$$\rho(x', \xi',z) =  \sqrt{z - r_0(x', \xi')} \in C^{\infty}(T^* \Gamma)$$
 be the root of the equation 
$$\rho^2 + r_0(x', \xi') - z = 0$$
  with $ \im \rho(x', \xi', z) > 0.$ For large $|\xi'|$ we have $|\rho(x', \xi', z)| \sim |\xi'|,\: \im \rho(x', \xi', z) \sim |\xi'|.$ Moreover, for $ z \in Z_1 \cup Z_3$  we have
$$\im \rho(x', \xi', z) \geq \frac{|\im z|}{2 |\rho|}, |\rho| \geq \sqrt{|\im z|},$$
while for $r_0 \geq 2$, we have
$$ C_1 \sqrt{r_0 + 1} \geq 2 \im \rho \geq |\rho| \geq C_2 \sqrt{r_0 + 1}.$$
For $z \in Z_2$ the last equality is true for all $(x', \xi')$ (see Lemma 3.1 in \cite{V}).

G. Vodev established for bounded domains $K \subset \R^d, \: d \geq 2,$ the following approximation of the interior Dirichlet-to-Neumann map $N_{int}(z, h)$ related to the boundary problem (\ref{eq:1.5}), where the equation $(-h^2 \Delta - z)u = 0$ is satisfied in $K$.
\begin{theorem} [\cite{V}, Theorem 3.3] For every $0 < \epsilon \ll 1$ there exists $0 < h_0(\epsilon) \ll 1$ such that for $z \in Z_1$ and $0 < h \leq h_0$ we have
\begin{equation} \label{eq:2.3}
\|\nc_{int}(z, h)(f) - Op_h(\rho + h b) f\|_{H_h^1(\Gamma)} \leq \frac{C h}{\sqrt{|\im z|}} \|f\|_{L^2(\Gamma)},
\end{equation}
where $b \in S^0_{0, 1}(\Gamma)$ does not depend on $h$ and $z$. Moreover, $(\ref{eq:2.3})$ holds for $z \in Z_2 \cup Z_3$ with $|\im z|$ replaced by $1.$
\end{theorem}

The same result remains true for unbounded domains $\Omega$ with $\nc_{int}(z, h)$ replaced by $\nc_{ext}(z, h)$ by modifications to the proof in \cite{V} based on the construction of a semi-classical parametrix close to the boundary. For reader convenience we recall below some facts from \cite{V} and we discuss some modifications which will be necessary for our exposition. Consider normal geodesic coordinates $(x_1, x')$ in a neighbourhood of a fixed point $x_0 \in \Gamma$, where $x_1 = \dist(x, \Gamma).$  
 Then $-h^2 \Delta- z$ in these coordinates has the form
$${\mathcal P}(z, h) = h^2D_{x_1}^2 + r(x,h D_{x'}) +  q(x, h D_x) + h^2 \tilde{q}(x)- z.$$
 with $D_{x_1} = -\ii \pa_{x_1},\: D_{x'} = -\ii \pa_{x'}, \: r(x, \xi') = \langle R(x)\xi', \xi'\rangle, \:q(x, \xi) = \langle q(x), \xi \rangle$.
Here $R(x)$ is a symmetric $(d-1) \times (d-1)$ matrix with smooth real-valued entries and $r(0, x', \xi') = r_0(x', \xi')$.
Let $\phi(\sigma) \in C^{\infty}(\R)$ be a cut-off function such that $\phi(\sigma) = 1$ for $|\sigma| \leq 1,$ $\phi(\sigma) = 0$ for $|\sigma| \geq 2.$
Let $\psi(x')$ be a $C^{\infty}$ cut-off function on $\Gamma$ supported in a small neighbourhood of $x_0$ and $\psi(x_0) = 1$. In \cite{V}, Proposition 3.4 for $\delta_1 > 0$ small enough one constructs a semi-classical parametrix
\begin{equation} \label{eq:2.4}
\begin{cases}\tilde{u}_{\psi}(x) = ( 2 \pi h)^{-d +1} \int\int e^{\frac{\ii}{h} \varphi(x, y', \xi', z)}\phi\Bigl(\frac{x_1}{\delta_1}\Bigr)\phi\Bigl(\frac{x_1}{\delta_1 \rho_1}\Bigr) a(x, \xi', z; h) f(y') dy'd\xi',\\
\tilde{u}_{\psi}\vert_{x_1 = 0} = \psi f,\end{cases}
\end{equation}
where $\rho_1 = 1$ if $ z \in Z_2 \cup Z_3,\: \rho_1 = |\rho|^3$ if $z \in Z_1$.  

 The phase $\varphi(x, y', \xi', z)$ is complex-valued and 
$$\varphi = - \langle x'- y', \xi'\rangle + \sum_{k=1}^{N-1} x_1^k \varphi_k(x', \xi', z) =  - \langle x'- y', \xi'\rangle + \tilde{\varphi},$$
$$a = \sum_{k=0}^{N-1} \sum_{j=0}^{N-1} x_1^k h^j a_{k,j}(x', \xi', z),$$
$N \gg 1$ being a large integer.

 Moreover, $\varphi_1 = \rho,\: \im \varphi \geq x_1 \im \rho/2,\: 0 \leq x_1 \leq 2\delta_1\min\{1, \rho_1\}$ and the amplitude $a$ satisfies $a\big\vert_{x_1 = 0} = \psi(x').$
 The phase $\varphi$ and the amplitude $a$ are determined so that
$$e^{-\frac{\ii\varphi}{h}} {\mathcal P}(z, h) e^{\frac{\ii\varphi}{h}}a = x_1^N A_N(x, \xi', z; h) + h^N B_N(x, \xi', z; h),$$
where $A_N,\: B_N$ are smooth functions. To describe the behavior of $A_N, B_N$, introduce the function $\chi(x', \xi') = \phi(\delta_0 r_0(x', \xi')),$ where $0 < \delta_0 \ll 1.$ Following \cite{V}, we say that a symbol $b \in C^{\infty}(T^*\Gamma)$ belongs to $S_{\delta_1, \delta_2}^{l_1}(\mu_1) + S_{\delta_3, \delta_4}^{l_2}(\mu_2)$ if
$$|\pa_{x'}^{\alpha} \pa_{\xi'}^{\beta} (\chi b)| \leq \cab |\mu_1|^{l_1 - \delta_1|\alpha|- \delta_2|\beta|},$$
$$|\pa_{x'}^{\alpha} \pa_{\xi'}^{\beta} ((1-\chi) b)| \leq \cab |\mu_2|^{l_2 - \delta_3|\alpha|- \delta_4|\beta|}, \: \forall \alpha, \forall \beta.$$
Therefore,
$$\pa_{x_1}^k A_N \in S_{2,2}^{2-3N-3k}(|\rho|) + S_{0, 1}^2(|\rho|),\: \pa_{x_1}^k B_N \in S_{2,2}^{3 - 4N-3k}(|\rho|) + S_{0,1}^{1-N}(|\rho|), \: \forall k \in \N$$
uniformly with respect to $z, h$ and $0 \leq x_1 \leq 2\delta_1\min\{1, \rho_1\}.$

For $z \in Z_{1, 0}$ and any integer $s \geq 0$, there exist $l_s, N_s> 0$ so that for $N \geq N_s$ we have the estimate (see Proposition 3.7 in \cite{V})
\begin{equation} \label{eq:2.5}
\|{\mathcal P}(z, h) \tilde{u}_{\psi}\|_{H_h^s(\Omega)} \leq C_N h^{-l_s} \Bigl(\frac{\sqrt{h}}{|\im z|}\Bigr)^{2N} \|f\|_{\lg},
\end{equation}
while for $z \in Z_2 \cup Z_3$ the above estimate holds with $|\im z|$ replaced by 1. Next introduce the operator
$$T_{\psi}(z, h) f: = D_{x_1} \tilde{u}_{\psi}\big\vert_{x_1 = 0} = Op_h(\tau_{\psi})f$$
with 
$$\tau_{\psi} = a \frac{\pa \varphi}{\pa x_1} \big\vert_{x_1 = 0} - \ii h \frac{\pa a}{\pa x_1}\big\vert_{x_1 = 0} = \psi \rho - \ii h \sum_{j=0}^{N-1} h^j a_{1, j}.$$
Let $G_D$ be the self-adjoin realization of the operator $-\Delta$ on $L^2(\Omega)$ with Dirichlet boundary condition on $\Gamma.$ Since the spectrum of $G_D$ is the positive real axis, for $z \in Z_1$ we have the estimate
$$\Bigl\|\Bigl(h^2 G_D - z)^{-1} \Bigr\|_{H_h^{2k} (\Omega) \to H_h^{2k}(\Omega)} \leq \frac{C_k}{|\im z|},\: \forall k \in \N,$$
while for $z \in Z_2 \cup Z_3$ the above estimate holds with $|\im z|$ replaced by 1. For $k = 0$ this estimate is trivial, and for $k \geq 1$ it follows from the coercive estimates for the Dirichlet problem in unbounded domains (see \cite{LP1})
$$\|v\|_{H_h^{2k}(\Omega)} \leq C_k'\Bigl( \|h^2 G_D v\|_{H_h^{2k - 2}(\Omega)} + \|v\|_{H_h^{2k-2}(\Omega)}\Bigr), \: v \in D(G_D) \cap H_h^{2k- 2}(\Omega).$$

Now let $u_{\psi} \in H_h^2(\Omega)$ be the solution of the problem 
$${\mathcal P}(z, h) u_{\psi} = 0\:\: {\text in}\:\Omega,\: u_{\psi}\vert_{\Gamma} = \psi f.$$
Then
$$w_{\psi}: = u_{\psi} - \tilde{u}_\psi + \Bigl(h^2 G_D - z\Bigr)^{-1} {\mathcal P}(z, h) \tilde{u}_{\psi}$$
will be a solution of $(h^2 G_D - z)w_{\psi} = 0$ in $\Omega$, $w_{\psi}\vert_{\Gamma} = 0$. Since for $z \in Z_1 \cup Z_2 \cup Z_3$ the point $z/h^2$ is not in the spectrum of $G_D$, one  deduces $w_{\psi} = 0$. This implies as in \cite{V} the following
\begin{proposition} For $z \in Z$ we have the estimate
\begin{equation} \label{eq:2.6}
\|\nc_{ext}(z, h)u_{\psi} - T_{\psi} (z, h) f\|_{H_h^1(\Gamma)} \leq C_N h^{-s_d} \Bigl(\frac{\sqrt{h}}{|\im z|}\Bigr)^{2N} \|f\|_{\lg}, \: \forall N \in \N 
\end{equation}
with constants $C_N, s_d > 0$, independent of $f, h$ and $z$, and $s_d$ independent of $N$. If $z \in Z_2 \cup Z_3$, then $(\ref{eq:2.6})$ holds with $|\im z|$ replaced by $1.$
\end{proposition}
Choose a partition of unity $\sum_{j = 1}^J \psi_j(x') = 1$ on $\Gamma$ and set $T(z, h) = \sum_{j=1}^J T_{\psi_j}(z, h)$. Notice that the principal symbol of $T(z, h)$ is $\rho.$
 By using Proposition 2.2 and repeating without any change the argument in Section 3 in \cite{V}, one concludes that the statement of Theorem 2.1 remains true replacing in (\ref{eq:2.3}) $\nc_{int}(z, h)$ by $\nc_{ext}(z, h)$.\\

\section{Eigenvalues-free regions in the case (A)}

In this section we suppose that $0 < \epsilon_0 \leq \gamma(x) \leq 1 - \epsilon_0,\: \epsilon_0 > 0, \: \forall x \in \Gamma.$ If $(u, v) \neq 0$ is an eigenfunction of $G$ with eigenvalue $\lambda \in \{z \in \C:\re z < 0\}$, then  $f = u\vert_{x \in \Gamma} \neq 0$. Indeed, if $f = 0$ on $\Gamma$, then $u \in H^2(\Omega)$ will be eigenfunction of the Dirichlet problem in $\Omega$ and this is impossible. From (\ref{eq:1.3}) one obtains the equation (\ref{eq:1.4}).

According to Theorem 2.1 with $\nc_{ext}(z, h)$, for $z \in Z_1,\: \delta = 1/2 - \ep,$ we have
\begin{equation} \label{eq:3.1}
\|Op_h(\rho)f  - \sqrt{z}\gamma f\|_{\lg} \leq C \frac{h}{\sqrt{|\im z|}} \|f\|_{\lg},
\end{equation}
where for $z \in Z_2 \cup Z_3$ the above estimate holds with $|\im z|$ replaced by 1. Here we use the fact that 
$$\|Op_h(b)\|_{L^2(\Gamma) \to L^2(\Omega)}  \leq C $$
which follows from \cite{V}, Proposition 2.1.
Introduce the symbol
$$c(x',\xi', z): = \rho(x',\xi', z)- \gamma\sqrt{z}.$$
We will show that $c(x', \xi', z)$ is elliptic in a suitable class. Write
$$ c(x', \xi', z) = \frac{(1 - \gamma^2) z - r_0(x', \xi')}{\rho(x', \xi', z) + \gamma \sqrt{z}}.$$

{\bf Case I.} $z \in Z_1.$
The symbol $c$ is elliptic for $|\xi'|$ large enough and  it remains to examine its behavior for $|\xi'| \leq C_0$. For these values of $\xi'$  we have $|\rho + \gamma \sqrt{z} | \leq C_1.$ 
First consider the set  
$${\mathcal F} = \{(x', \xi'):\:
 |1- r_0(x'. \xi')| \leq \frac{\epsilon_0^2}{2}\}.$$
Then $\re \Bigl((1 - \gamma^2) z - r_0 \Bigr) = 1- r_0 - \gamma^2 \leq -\frac{\epsilon_0^2}{2}.$
If $(x', \xi') \notin {\mathcal F}$, we get
$$\im \Bigl( (1 - \gamma^2) z - r_0 \Bigr) = (1 - \gamma^2) \im z  \geq (1 - \gamma^2) h^\delta \geq \ep_0 h^\delta.$$
Consequently, the symbol $c$ is elliptic and
$$\im (\rho + \gamma \sqrt{z}) = \im \rho + \gamma \im \sqrt{z} \geq C h^{\delta}.$$ 
Hence, for bounded $|\xi'|$ we have $|c| \geq C_3 h^{\delta}, C_3> 0,$
while for large $|\xi'|$ we have $|c| \sim |\xi'|.$ As in Section 2 we use the function $\chi$ 
and define ${\mathcal M}_1: = Z_1 \times \supp \chi,\: {\mathcal M}_2: = (Z_1 \times \supp (1 - \chi)) \cup ((Z_2\cup Z_3) \times T^* \Gamma).$
Set $\langle \xi'\rangle = ( 1 + |\xi'|)^{1/2}.$ It is easy to see that for $(z, x', \xi') \in {\mathcal M}_1,$ we have
\begin{equation} \label{eq:3.2}
\big | \pa_{x'} ^{\alpha} \pa_{\xi'}^{\beta} \rho \big | \leq C_{\alpha, \beta} |\im z|^{1/2 - |\alpha| - |\beta|},\: |\alpha| + |\beta| \geq 1,
\end{equation}
$|\rho| \leq C,$ while for $(z, x', \xi') \in {\mathcal M}_2$ we have
\begin{equation} \label{eq:3.3}
\big | \pa_{x'} ^{\alpha} \pa_{\xi'}^{\beta} \rho \big | \leq C_{\alpha, \beta} \langle \xi' \rangle ^{1- |\beta|}.
\end{equation}
Thus, we conclude that $c =(\rho - \gamma \sqrt{z})\in S^{0, 1}_{\delta}.$\\

Now consider the symbol $c^{-1} = \frac{\rho + \gamma \sqrt{z}}{(1 - \gamma^2) z - r_0}.$ Since $\rho + \gamma \sqrt{z} \in S^{0, 1}_{\delta},$ it remains to study the properties of $g := ((1- \gamma^2)z - r_0)^{-1}.$  For $(x', \xi') \in {\mathcal F}$, we get $|(1 - \gamma^2) z - r_0| \geq \frac{\epsilon_0^2}{2} > 0$ and 
$$|\pa_{x'}^{\alpha} \pa_{\xi'}^{\beta} g| \leq \cab .$$
Therefore for $(x', \xi') \in {\mathcal F}$, we have
\begin{equation} \label{eq:3.4}
|\pa_{x'}^{\alpha} \pa_{\xi'}^{\beta} (c^{-1})| \leq \cab |\im z|^{1/2 - |\alpha| - |\beta|}.
\end{equation}

 Next for $(x', \xi') \notin {\mathcal F}$ notice   that for every $0 < \delta' \ll 1,$ if $|( 1 - \gamma^2)  - r_0| \leq \delta',\: \im z \neq 0,$  we have
\begin{equation} \label{eq:3.5}
\big | \pa_{x'} ^{\alpha} \pa_{\xi'}^{\beta} g| \leq C_{\alpha, \beta} |\im z|^{-1 - |\alpha| - |\beta|},
\end{equation}
while for $|( 1 - \gamma^2)  - r_0| \geq \delta'$ we get
\begin{equation} \label{eq:3.6}
\big | \pa_{x'} ^{\alpha} \pa_{\xi'}^{\beta} g | \leq C_{\alpha, \beta}\langle \xi'\rangle ^{-2 - |\beta|}.
\end{equation}
On the other hand, $(x', \xi') \notin {\mathcal F}$ yields $|1 - r_0(x', \xi')| > \frac{\epsilon_0^2}{2}$ and for $(x', \xi') \notin {\mathcal F}$ we
obtain
\begin{equation} \label{eq:3.7}
|\pa_{x'}^{\alpha} \pa_{\xi'}^{\beta} \rho| \leq \cab \langle \xi' \rangle ^{1- |\alpha| - |\beta|}.
\end{equation}
Thus for bounded $|\xi'|$ and $(x',\xi') \notin {\mathcal F}$, we deduce
\begin{equation} \label{eq:3.8}
|\pa_{x'}^{\alpha} \pa_{\xi'}^{\beta}(c^{-1})| \leq \cab |\im z|^{-1- |\alpha| - |\beta|}.
\end{equation}
Combining this with the estimates (\ref{eq:3.4}), one concludes that $|\im z| c^{-1}\in S_{\delta}^{0, -1}.$ \\

{\bf Case II.} $z \in Z_2$. We have
$$ \re \Bigl((1- \gamma^2)  z - r_0\Bigr) \leq -(1 - \gamma^2) \leq  -\ep_1 < 0.$$
Consequently,  $c$ is elliptic and  $c \in S_0^{0, 1}, c^{-1} \in S_0^{0, -1}.$ \\

{\bf Case III.} $z \in Z_3$.  In this case $\im z = 1$ and one has
$$\big|\im \Bigl((1-\gamma^2)z - r_0\Bigr)\big| = |(1- \gamma^2)| \geq \epsilon_1 > 0.$$ 
This implies that $c \in S_0^{0, 1}$ is elliptic and  $c^{-1} \in S_0^{0, -1}.$\\

Consequently, we get
$$\|Op_h(c^{-1}) g \|_{L^2(\Gamma)} \leq  C |\im z|^{-1} \|g\|_{L^2(\Gamma)}$$
and,  applying  (\ref{eq:3.1}), we deduce
$$\|Op_h(c^{-1}) Op_h(c) f \|_{L^2(\Gamma)} \leq C_5 \frac{h}{|\im z|^{3/2}} \|f\|_{L^2(\Gamma)}. $$
On the other hand, for $|\alpha_1| + |\beta_1| \geq 1,\: |\alpha_2| + |\beta_2| \geq 1$ and $|\xi'| \leq C_0$ according to (\ref{eq:3.2}), (\ref{eq:3.4}), (\ref{eq:3.7}), (\ref{eq:3.8}), for $(x', \xi') \in {\mathcal F}$ we get
\begin{equation} \label{eq:3.9}
\Bigl|\pa_{x'}^{\alpha_1} \pa_{\xi'}^{\beta_1} (c^{-1}(x', \xi')) \partial_{x'}^{\alpha_2} \partial_{\xi'}^{\beta_2} c(x', \xi')\Bigl|\\
\leq  C_{\alpha_1, \beta_1, \alpha_2, \beta_2} |\im z|^{1 -(|\alpha_1| + |\beta_1| + |\alpha_2| + |\beta_2|)},
\end{equation}
while for $(x', \xi') \notin {\mathcal F}$ we have
\begin{equation} \label{eq:3.10}
\Bigl|\pa_{x'}^{\alpha_1} \pa_{\xi'}^{\beta_1} (c^{-1}(x', \xi')) \pa_{x'}^{\alpha_2} \pa_{\xi'}^{\beta_2} c(x', \xi')\Bigl|
\leq C_{\alpha_1, \beta_1, \alpha_2, \beta_2} |\im z|^{-1 - (|\alpha_1| + |\beta_1|)}.
\end{equation}
Consider the operator $Op_h(c^{-1})Op_h(c) - I$. Following Section 7 in \cite{DS}, the symbol of this operator is given by
$$\sum_{j=1}^{N} \frac{(i h)^j}{j !} \sum_{|\alpha| = j} D_{\xi'}^{\alpha} (c^{-1})(x', \xi') D_{y'}^{\alpha}c(y', \eta')\big\vert_{x'= y', \xi'= \eta'} + \tilde{b}_{N}(x', \xi')$$
$$= b_N(x', \xi') + \tilde{b}_N(x', \xi'),$$
where
$$|\pa_{x'}^{\alpha} \tilde{b}_N(x', \xi')| \leq C_{\alpha} h^{N(1 - 2 \delta) -s_d- |\alpha|/2}.$$
Applying (\ref{eq:2.2}), one deduces for $N$ large enough
$$\|Op_h(\tilde{b}_N)\|_{\lg \to \lg} \leq C h.$$

On the other hand, the estimates (\ref{eq:3.9}), (\ref{eq:3.10}) yield
$$|\pa_{x'}^{\alpha} b_N(x', \xi')| \leq C_{\alpha} \frac{h}{|\im z|^2}h^{-\delta|\alpha|}.$$
Thus, applying once more (\ref{eq:2.2}), one gets
$$\|Op_h(c^{-1}) Op_h(c) f - f \|_{L^2(\Gamma)} \leq C_6 \frac{h}{|\im z|^{2}} \|f\|_{L^2(\Gamma)}.$$
A combination of the above estimates implies
\begin{equation} \label{eq:3.11}
\|f\|_{L^2(\Gamma)} \leq C_7 \Bigl(  h^{1-2\delta} + h^{1-\frac{3}{2}\delta}\Bigr) \|f\|_{L^2(\Gamma)}.
\end{equation}
Since $\delta = 1/2 - \epsilon,\: 0 < \epsilon \ll 1$, for $0 < h \leq h_0(\epsilon)$ small enough (\ref{eq:3.8}) yields $f = 0.$\\

 Going back to $\lambda = \frac{\ii \sqrt{z}}{h},$ we have
$$\re \lambda = - \frac{\im \sqrt{z}}{h}, \: \im \lambda = \frac{\re \sqrt{z}}{h}.$$
Suppose that $z \in Z_1.$ Then

$$|\re \lambda| \geq C (h^{-1})^{1- \delta},\: |\im \lambda| \leq C_1 h^{-1} \leq C_2 |\re \lambda|^{\frac{1}{1-\delta}}.$$
So if $|\re \lambda| \geq C_3 |\im \lambda|^{1 - \delta}, \: \re \lambda \leq -C_4 < 0$ there are no eigenvalues $\lambda = \frac{\ii \sqrt{z}}{h}$ of $G$.
For $z \in Z_2 \cup Z_3$ there are no eigenvalues $\lambda$ too if $|\lambda| \geq R_0.$

This shows that in the case (A) for every $0 < \epsilon \ll 1$  the eigenvalues of $G$ must lie in the region $\Lambda_{\epsilon}$ defined in Theorem 1.2.

\section{Eigenvalues-free region in the case (B)}
\def\lg{L^2(\Gamma)}
In this section we deal with the case (B). 
The  analysis of Section 3  works only for $z \in Z_1 \cup Z_3$.
Indeed for $z \in Z_1$  we have 
$$ \re ((1 - \gamma^2)- r_0) \leq (1 -\gamma^2) < - \eta_0 < 0.$$ 
 The symbol $g$ introduced in the previous section satisfies the estimates (\ref{eq:3.5}) and   $c \in S_{\delta}^{0, 1}, c^{-1} \in S_{\delta}^{0, -1}.$ For $z \in Z_3$ we apply the same argument. Thus for $z \in Z_1 \cup Z_3$ we obtain that the eigenvalues $\lambda = \frac{\ii \sqrt{z}}{h}$of $G$  must lie in $\Lambda_{\epsilon}.$
For $z \in Z_2$ the argument exploited in the case (A) breaks down since for $\re z = -1, \im z = 0$ the symbol
$$\ii [1+ r_0(x',\xi') - \gamma(x')]$$
is not elliptic and it may vanish for some $(x_0', \xi_0').$

 In the following we  suppose that $z \in Z_2$. Therefore Proposition 2.2 yields a better approximation
\begin{equation} \label{eq:4.1}
\| \nc_{ext}(z, h)(f) - T(z, h) f\|_{H^{1}(\Gamma)} \leq C_N h^{-s_d + N} \|f\|_{L^2(\Gamma},\: \forall N \in \N. 
\end{equation}
 If $f \neq 0$ is the trace of an eigenfunction of $G$, from the equality (\ref{eq:1.4}) we obtain
$$|\re \langle T(z, h) f - \gamma \sqrt{z} f, f \rangle_{L^2(\Gamma)} | \leq C_N h^{-s_d + N} \|f\|_{\lg}.$$
There exists $t$ with $0 < t < 1$ such that
\begin{eqnarray} \label{eq:4.2}
\re \langle \Bigl(T(z, h) - \gamma \sqrt{z}\Bigr)f, f \rangle_{\lg}  = \re \langle T(-1, h)f, f \rangle_{\lg} \nonumber  \\
-  \im z \im \Bigl[ \langle \frac{\partial T}{\partial z} (z_t, h) - \gamma \frac{1}{2 \sqrt{z_t}}f, f \rangle_{\lg} \Bigr]
\end{eqnarray}
with $z_t = - 1 + \ii t \im z \in Z_2.$ The next Lemma is an analogue of Lemma 3.9 in \cite{V}.
\begin{lemma} Let $z \in Z_2$ and let $f = u\vert_{\Gamma}$ be the trace of an eigenfunction $u$ of $G$ with eigenvalue $\lambda = \frac{\ii \sqrt{z}}{h}.$ Then
\begin{equation} \label{eq:4.3}
\Bigl\|\frac{d T}{dz} (z, h) f - Op_h\Bigl( \frac{d\rho}{dz}(z)\Bigl)f \Bigr\|_{\lg} \leq C h \|f\|_{H_h^{-1}(\Gamma)}
\end{equation}
with a constant $C > 0$ independent of $z, h$ and $f$. Moreover,
\begin{equation} \label{eq:4.4}
|\re \langle T(-1, h)f , f \rangle_{\lg}| \leq C_N h^{-s_d + N} \|f\|_{\lg}, \: \forall N \in \N.
\end{equation}
\end{lemma}
 {\it Proof.} The proof of (\ref{eq:4.3}) is the same as in \cite{V} since for $z \in Z_2$ we get
$$\sum_{j=0}^{N-1} h^j \frac{d a_{1, j}}{dz} \in S_0^{0, -1}.$$

 To establish (\ref{eq:4.4}), we apply Green's formula in the unbounded domain $\Omega$.
By using the notation of Section 3, set $\tilde{u} = \sum_{j=1}^J \tilde{u}_{\psi}.$ Then $-\ii h  \pa_{\nu}\tu\vert_{\Gamma} = T(z, h) f$ and for $R \gg 1$ the function $\tu$ vanishes for $|x| \geq R$. Thus one obtains
$$\ii \langle\Delta \tu, \tu \rangle_{L^2(\Omega)} = - \ii\int_{\Omega} |\nabla \tu|^2 dx  - \ii\langle \pa_{\nu}\tu, \tu \rangle_{\lg}.$$
 Multiplying the above equality by $h$ and taking  the real part, we deduce
$$- \im h \langle \Delta \tu, \tu \rangle_{L^2(\Omega)} = \re \langle T(z, h)f, f \rangle_{\lg}.$$
Therefore,
$$\re \langle T(z, h)f, f \rangle_{\lg}. = -\im h \langle (\Delta - h^{-2})\tu, \tu \rangle_{L^2(\Omega)}= -\im h^{-1} \langle {\mathcal P}(-1, h) \tu, \tu \rangle_{L^2(\Omega)}$$
and
$$|\re \langle T(z, h)f, f \rangle_{\lg}| \leq h^{-1} \|{\mathcal P}(-1, h) \tu\|_{L^2(\Omega)} \|\tu\|_{L^2(\Omega)}.$$
It is easy to see that $\|\tu\|_{L^2(\Omega)} \leq C h^{-s_d}\|f\|_{L^2(\Omega)}$ and combining this with (\ref{eq:2.5}) for $z \in Z_2$, we obtain (\ref{eq:4.4}).  
\hfill\qed\\

From (\ref{eq:4.2}), (\ref{eq:4.4}) and $\im z \neq 0$ we have
\begin{equation} \label{eq:4.5}
| \im \langle \Bigl(\frac{\partial T}{\partial z} (z_t, h) - \frac{\gamma}{2 \sqrt{z_t}}\Bigr) f, f \rangle_{\lg}| \leq C_N \frac{h^{-s_d + N}}{|\im z|} \|f\|_{\lg}.
\end{equation}
Consider the operator $L := Op_h(\frac{d\rho}{dz}(z_t)) -  \frac{\gamma}{2 \sqrt{z_t}}$ and notice that
\begin{equation} \label{eq:4.6}
\Bigl|\im \langle \Bigl(\frac{\partial T}{\partial z} (z_t, h) - \frac{\gamma}{2 \sqrt{z_t}}\Bigr) f, f \rangle_{\lg} 
- \im \langle Lf, f \rangle_{\lg} \Bigr| \leq C h \|f\|_{\lg}.
\end{equation}
On the other hand,
$$\im \langle Lf, f \rangle _{\lg} = \frac{1}{ 2\ii} \langle (L - L^*)f, f \rangle _{\lg}$$
and the principal symbol of $\frac{1}{2 \ii} (L- L^*)$ has the form
$$ s(x', \xi'; z) := \frac{1}{2}\im \Bigl[ \frac{1}{\sqrt{ - 1+ \ii t \im z - r_0}} - \frac{\gamma} {\sqrt{-1 + \ii t \im z}}\Bigr].$$
Let 
$$z_t = y e^{\ii (\pi - \varphi)},\:y = \sqrt{1 + t^2 |\im z|^2},\: |\varphi| \leq \pi/4.$$
Here and below we omit the dependence of $y$ on $t$.
 Then $ 1 \leq y \leq \sqrt{2}$ and 
$$\sqrt{z_t} = \sqrt{y} \sin \varphi/2 + \ii \sqrt{y} \cos \varphi/2,\: \im \frac{1}{\sqrt{z_t}} = - \frac{ \cos\varphi/2}{\sqrt{y}}.$$
In the same setting
$$-1 + \ii t \im z - r_0 = q e^{\ii(\pi - \psi)},\:q = \sqrt{(1 + r_0)^2 + t^2(\im z)^2}, |\psi| \leq \pi/4,$$
we see that
$$\im \frac{1}{\sqrt{ - 1+ \ii t \im z - r_0}} = -\frac{\cos \psi/2}{\sqrt{q}}.$$
Therefore
$$s = \frac{1}{2\sqrt{y q}} \Bigl( \gamma \sqrt{q} \cos\varphi/2 - \sqrt{y} \cos \psi/2\Bigr)= \frac{ \gamma^2 q \cos^2 \varphi/2 - y \cos^2 \psi/2}{ 2\sqrt{yq} ( \gamma \sqrt{q} \cos \varphi/2 + \sqrt{y} \cos \psi/2)}.$$
To prove that $s$ is elliptic, it is sufficient to show that 
$$ \gamma^2 q( 1+ \cos \varphi) - y(1 + \cos \psi)= \gamma^2 q\Bigl( 1+ \frac{1}{y}\Bigr) - y\Bigl(1 + \frac{1 + r_0}{q}\Bigr)$$
$$= \frac{1}{y q} \Bigl[ \gamma^2 q^2( 1 + y) - y^2( 1 + q + r_0)\Bigr]$$
is elliptic. Consider the function
$$F(r_0) = \gamma^2\Bigl( (1 + r_0)^2 + t^2\im^2 z\Bigr)( 1 + y) - y^2 \Bigl( 1 + \sqrt{( 1 + r_0)^2 + t^2\im^2 z} + r_0\Bigr).$$
Clearly,
$$F(0) = ( \gamma^2 - 1)( 1 + y) y^2 \geq \eta_1 > 0,$$
since in the case (B) we have $\gamma^2 - 1 \geq \eta_0 > 0.$ Next, for $\gamma \geq 1, \:r_0 \geq 0$ we have
$$\frac{\partial F}{\partial r_0} = 2 \gamma^2 ( 1+ y) (1 + r_0) - y^2\Bigl( 1 + \frac{1 + r_0}{\sqrt{(1+r_0)^2 + t^2\im^2 z}}\Bigr)$$
$$ \geq 2 \Bigl(\gamma^2 ( 1+ y) (1 + r_0) -  y^2\Bigr) \geq 2( 1 + y - y^2).$$
On the other hand, it is clear that $1 +y - y^2 > 0$ for $0 \leq y < \frac{ 1 + \sqrt{5}}{2}.$ In our case $1 \leq y \leq \sqrt{2} < \frac{ 1 + \sqrt{5}}{2}$ and we deduce $\frac{\partial F}{\partial r_0}(r_0) > 0$ for $r_0 \geq 0,\: 1 \leq y \leq \sqrt{2}.$ This implies $F(r_0) > 0$ for $r_0 \geq 0$ and $s$ is elliptic. Consequently, 
$$ \im \langle Lf, f \rangle_{\lg} \geq (\eta_2 - C h) \|f\|_{\lg},\: \eta_2 > 0$$
and for small $h$ and $\|f\|_{\lg} \neq 0,\: \im z \neq 0,$ we deduce from (\ref{eq:4.5}) and (\ref{eq:4.6})
$$|\im z| \leq C_N' h^{-s_d + N}\leq B_N h^N,\: \forall N \in \N.$$
Going back to $\lambda = \frac{\ii \sqrt{z}}{h},$ we have
$$\re \sqrt{z} = \mu^{1/2} \sin \varphi/2, \: \im \sqrt{z} = \mu^{1/2} \cos \varphi/2,\: \mu = \sqrt{1 + (\im z)^2} \leq \sqrt{2}$$
and $0 \leq \sin \varphi \leq B_N h^{N}.$ This implies for $h$ small enough the estimate
$$|\im \lambda| = \Bigl| \frac{ \re \sqrt{z}}{h} \Bigr| \leq 2^{1/4}B_N (h^{-1})^{-N + 1}\leq C_N |\re \lambda|^{-N}.$$
Thus for $z \in Z_2$ and every $N \in \N$ the eigenvalues $\lambda = \frac{\ii \sqrt{z}}{h}$ of $G$ lie in ${\mathcal R}_{N}$ and this completes the proof of Theorem 1.2.\\
\hfill\qed

The eigenvalues of $G$ could have accumulation points on $\ii \R$. For odd dimension $d$ Lax and Phillips \cite{LP} proved that the scattering matrix $S(z)$ is invertible for $z = 0.$ This leads easily to the following

\begin{proposition} Assume $d$ odd. The operator $G$ has no a sequence of eigenvalues $\lambda_j,\: \re \lambda_j < 0$ such that $\lim_{j\to \infty} \lambda_j= \ii z_0,\: z_0 \in \R.$
\end{proposition}
The proof is the same as that of Proposition 4.11 in \cite{CPR}. The above proposition does not exclude the possibility to have eigenvalues $\lambda_j$ with $|\im \lambda_j| \to +\infty.$ On the other hand, Theorem 1.3, established in the next section, implies that for strictly convex obstacles and $\gamma(x) > 1$ the imaginary part of all eigenvalues of $G$ is bounded by a constant $R_0 > 0$ and for $d$ odd we can apply Proposition 4.2.

\section{Eigenvalue-free region for strictly convex obstacles in the case (B)}

In this section we study the eigenvalues-free regions when $K$ is a strictly convex obstacle. Let $0 < \ep \ll 1/2$ be a small number. 
Set
$$\chi_1(x', \xi') = \phi \Bigl(\frac{1 - r_0(x', \xi')}{h^{\ep/2}}\Bigr),$$
where $\phi$ is the function introduced in Section 2. Notice that on the support of $1 - \chi_1$ we have $|1 - r_0(x', \xi')| \geq h^{\ep/2}.$
By a modification of the construction in \cite{V} (see also \cite{PV}) we can construct a semi-classical parametrix $\tu_{\psi}$ having the form (\ref{eq:2.4}), where $\psi f$ is replaced by $)p_h(1- \chi_1) \psi f.$

Then for $|1 - r_0(x', \xi')| \geq h^{\ep/2}$  we have $|\rho |^2 \geq h^{\ep/2}$ and we can improve the estimate (\ref{eq:2.5}) obtaining
\begin{equation} \label{eq:4.7}
\|{\mathcal P}(z, h) \tilde{u}_{\psi}\|_{H^s(\Omega)} \leq C_N h^{-l_s} \Bigl(\frac{h}{h^{\ep/2}|\im z|}\Bigr)^{N} \|f\|_{\lg}, \: |\im z| \geq h^{1- \ep}.
\end{equation}
To do this, one repeats without changes the argument in Section 3 of \cite{V} replacing the lower bound $|\rho|^2 \geq |\im z|$ by $|\rho|^2 \geq h^{\ep/2}.$ Consequently, the right hand side of (\ref{eq:4.7}) is estimated by ${\mathcal O}_{N}(h^{-l_s + N\ep/2})$ and this yields a semi-classical parametrix
$$ {\mathcal P}(z, h) w_1 ={\mathcal O}_N(h^N),\: w_1\vert_{x \in \Gamma} = Op_h(1 - \chi_1) \psi f.$$
Consider a partition of unity $\chi_{\delta}^{-} + \chi_{\delta}^0 + \chi_{\delta}^{+} = 1$ on $T^*(\Gamma)$, where the functions $\chi_{\delta}^-, \chi_{\delta}^0, \chi_{\delta}^+ \in S_{\delta, 0}^0$ are with values in $\R^+$ and such that supp $\chi_{\delta}^- \subset \{r_0 - 1\leq -h^{\delta}\},$
supp $\chi_{\delta}^+ \subset  \{r_0 - 1 \geq h^{\delta}\},$ supp $\chi_{\delta}^0 \subset \{|r_0 - 1| \leq 2 h^{\delta}\},\: \chi_{\delta}^0 = 1$ on $\{|r_0 -1| \leq h^{\delta}\}.$ Then, as in \cite{V}, \cite{V1}, we obtain the following
\begin{theorem} [Theorem 2.1, \cite{V1}] For every $0 < \epsilon \ll 1$ there exists $h_0(\epsilon) > 0$ such that for $0 < h \leq h_0(\epsilon), |\im z| \geq h^{1- \epsilon},$ we have
\begin{equation}
\|N_{ext}(z, h) Op_h(\chi_{\ep/2}^{-}) - Op_h(\rho \chi_{\ep/2}^{-})\|_{\lg \to \lg} \leq C h^{1/2}
\end{equation} \label{eq:5.2}
and for $|\im z| \leq h^{\ep}$ we have the estimate
\begin{equation} \label{eq:5.3}
\|N_{ext}(z, h) Op_h(\chi_{\ep/2}^{+}) - Op_h(\rho \chi_{\ep/2}^{+})\|_{\lg \to \lg} \leq C h^{1/2}.
\end{equation}
\end{theorem}

Thus the problem is to get an estimate for $\|N_{ext}(z, h) Op_h(\chi_{\ep/2}^0)\|_{\lg \to \lg}.$ We will prove the following
\begin{theorem}  For $h^{2/3} \leq \im z \leq h^{\ep}$ we have the estimate
\begin{equation} \label{eq:5.4}
\|\nc_{ext}(z, h) Op_h(\chi_{\ep/2}^0)\|_{\lg \to \lg} \leq C h^{\ep/2}.
\end{equation}
\end{theorem}
\begin{remark} By the analysis in \cite{V1} we may cover the region $h^{1-\ep} \leq \im z \leq h^{\ep}$ but the above result is sufficient for our analysis since the region $0 < \im z \leq h^{2/3}$ is examined in Chapter 9 and 10 in \cite{Sj}, where a parametrix for the exterior Dirichlet problem is  constructed with a precise estimate of the symbol of $\nc_{ext}$ in small neighbourhood of the glancing set (see $(10.31)$ in \cite{Sj}). 
\end{remark}
Set for simplicity of notation $\mu = \im z.$ We will follow closely the construction of a semi-classical parametrix  in Sections 5, 6 in \cite{V1}. The only difference is that we deal with an unbounded domain and the local form of ${\mathcal P}$ slightly changes. For the convenience of the reader we are going to recall the result in \cite{PoV}. Let $\Omega_{\delta} =  \{x \in \Omega:\: \dist(x, \Gamma) < \delta\}$. Since $K$ is strictly convex, in local normal geodesic coordinates $(x, \xi) \in T^*(\Omega_{\delta})$, considered in Section 2, the principal symbol of ${\mathcal P}$ becomes
$$p(x, \xi) = \xi_1^2 + r_0(x', \xi') + x_1 q_1(x', \xi') -1 - \ii \mu+ \oc( x_1^2 r_0)$$
with $0 < C_1 \leq q_1(x', \xi') \leq C_2$. Here locally in the interior of $K$ we have $x_1 > 0,$ while in the exterior of $K$ we have
$x_1 < 0.$
Following \cite{V1}, denote by ${\mathcal R}$ the set of functions $ a \in C^{\infty}(T^*(\Omega_{\delta}))$ satisfying with all derivatives the
estimates
$$a = \oc (x_1^{\infty}) + \oc (\xi_1^{\infty}) + \oc( ( 1- r_0)^{\infty})$$
in a neighbourhood of ${\mathcal K}: = \{(x, \xi): x_1 = \xi_1 = 1 - r_0 = 0\}.$ It was shown in Theorem 3.1 in \cite{PoV} that there exists an exact symplectic map $\chi: T^* (\Omega_{\delta}) \to T^* (\Omega_{\delta})$ so that $\chi(x, \xi) = (y(x, \xi), \eta(x, \xi))$ satisfies
$$y_1 = x_1 q_1(x', \xi')^{-1/3} + \oc(x_1^2) + \oc(x_1(1- r_0)),$$
$$\eta_1 = \xi_1 q_1(x', \xi')^{1/3} + \oc(x_1) + \oc(\xi_1(1- r_0)),$$
$$(y', \eta') = (x', \xi') + \oc(x_1),$$
$$(p \circ \chi(x, \xi)) = \Bigl(q_1(x', \xi')^{2/3} + \oc(x_1)\Bigr)(\xi_1^2 + x_1 - \zeta(x', \xi')),\: ({\rm mod}\:{\mathcal R})$$
in a neighbourhood of ${\mathcal K}$ with
$$ \zeta(x', \xi') = \Bigl(q_1(x', \xi')^{-2/3} + \oc(1- r_0)\Bigr) (1 + \ii \mu - r_0(x', \xi')).$$
 Let ${\mathcal U} \subset T^*(\Omega_{\delta})$ be a small neighbourhood of ${\mathcal K}.$ By using a $h$-Fourier integral operator on $\Omega_{\delta}$ associated to the canonical relation
$$\Lambda = \{ (y, \eta, x, \xi) \in T^* (\Omega_{\delta}) \times T^* (\Omega_{\delta}):\: (y, \eta) = \chi(x, \xi),\: (x, \xi) \in {\mathcal U}\},$$
one transforms ${\mathcal P}$ into an operator $P_0'$ which in the new coordinates denoted again by $(x, \xi)$ has the form
$$P_0'= D_{x_1}^2 + x_1 - L_1(x', D_{x'}; h) - \ii \mu L_2(x', D_{x'}; h),$$
where $L_j(x', \xi'; h) = \sum_{k=0}^{\infty} h^k L_j^{(k)}(x', \xi'),\: j = 1, 2,$ with 
$$L_1^{(0)}(x', \xi') =  \Bigl(q_1(x', \xi')^{-2/3} + \oc(1- r_0)\Bigr)(1- r_0(x', \xi')),$$
$$L_2^{(0)}(x', \xi') = q_1(x', \xi')^{-2/3} + \oc(1 - r_0).$$
By a simple change of variable $t = - x_1$, we pass to the situation when the exterior of $K$ is presented by $t > 0$. Next one applies a new symplectic transformation of the tangential variables $(x^\#, \xi^\#) = \chi^{\#}(x', \xi') \in T^*(\Gamma)$ so that $\xi_d^{\#} = - L_1^{(0)}(x', \xi')$ (see Section 2 in \cite{V1}). Therefore the operator $P_0'$ is transformed into
\begin{equation} \label{eq:5.5}
\tilde{P}_0 = D_{t}^2 -t + D_{x_d^\#} - \ii \mu q(x^\#, D_{x^\#}) + {\mathcal Q}(x^\#, D_{x^\#}; \mu, h),
\end{equation}
where $q(x^\#, \xi^\#) > 0,\: q \in S_0^0$ in a neighbourhood of $\xi_d^\# = 0$ and 
$${\mathcal Q} = \sum_{k=1}^{\infty} h^k Q_k(x^\#, \xi^\#;\mu).$$ 
 The only difference with \cite{V1} is the sign $(-)$ in front of $t$ in the form of $\tilde{P}_0.$\\
For simplicity of the notations we denote the coordinates $(x^\#, \xi^\#)$ by $(y, \eta)$ and consider the operator
$$P_0 = D_t^2 - t + D_{y_d} - \ii \mu q(y, D_y) + h \tilde{q}(y, D_y; \mu, h)$$
with $0 < C_1 \leq q(y, \eta) \leq C_2,\: q \in S_0^0,\: \tilde{q} \in S_0^0.$ Notice that we have the term $-\ii \mu q(y, \eta)$ with $\mu > 0$, while in \cite{V1} the model operator involves $\ii\mu q(y, \eta)$ since the sign of $\mu$ is not important for the argument in Sections 5, 6 of \cite{V1}.

First we will treat  the situation examined in Section 6 in \cite{V1}
when $\mu > 0$ and $\eta_d$ satisfy the conditions
\begin{equation} \label{eq:5.6} 
\mu \sqrt{\mu  + |\eta_d|} \geq h^{1- \ep},
\end{equation}
\begin{equation} \label{eq:5.7}
\mu + |\eta_d| \leq \oc(h^{\ep}).
\end{equation}
Clearly, if $h^{2/3 - \ep} \leq \mu \leq h^{\ep},$ the condition (\ref{eq:5.6}) holds. The same is true also if $h^{2/3} \leq \mu \leq h^{2/3 -\ep}$ and $|\eta_d| \geq h^{1/2 - \ep}.$ Introduce the function 
$$\Phi_1(\eta_d) = \begin{cases} \phi(\frac{\eta_d}{h^{\ep}}), \:{\rm if}\: \mu \geq h^{2/3 - \ep},\\
\Bigl(1 - \phi(\frac{\eta_d}{h^{1/2 - \ep}})\Bigr) \phi(\frac{\eta_d}{h^{\ep}}),\:{\rm if}\: h^{2/3} \leq \mu < h^{2/3- \ep}, \end{cases}$$
where $\phi$ is the function introduced in Section 2.
 Let $\rho$ be the solution of the equation
$$\rho^2 + \eta_d - \ii\mu q(y, \eta) = 0$$
with $\im \rho > 0.$ With a minor modifications of the argument in Section 6 in \cite{V1} we may construct a parametrix 
$\tilde{u}_1 = Op_h(A(t))f$, where
$$A(t) = \phi\Bigl(\frac{t}{\delta_1 |\rho|^2}\Bigr) a(t, y, \eta; \mu, h) e^{\frac{\ii \varphi(t, y, \eta: \mu)}{h}}$$ 
and $\delta_1 > 0$ is small enough.  We take $\varphi$ and $a$ in the form
$$\varphi = \sum_{k=1}^M t^k \varphi_k,\: a = \sum_{0 \leq k + \nu \leq M} h^k t^{\nu} a_{k, \nu},$$
where $M \gg 1$ and $\varphi_k$ and $a_{k, \nu}$ do not depend on $t$. We choose $a_{0,0} = \Phi_1(\eta_d),\: a_{k, 0} = 0$ for $k \geq 1$. We have the identity
$$e^{-\ii \varphi/h}( D_t^2 - t + \eta_d - \mu q(y, \eta)- \ii h \pa_{y_d}) (e^{\ii \varphi/h} a)$$
$$= - 2 \ii h \pa_t \varphi \pa_t a - h^2 \pa_t^2 a - \ii h \pa_{y_d} a + ((\pa_t \varphi)^2 + \pa_{y_d}\varphi - t - \rho^2)a$$
$$ = -2 \ii h \sum_{0 \leq k + \nu \leq 2M - 2} h^k t^{\nu} \sum_{j = 0}^{\nu} (j+1)(\nu + 1 - j) \varphi_{\nu + 1 - j} a_{k, j+1}$$
$$ - h \sum_{0 \leq k + \nu \leq M - 1} (\nu + 1)(\nu + 2) h^k t^{\nu} a_{k-1, \nu + 2} - \ii h \sum_{0 \leq k + \nu \leq M} h^k t^{\nu} \pa_{y_d}
a_{k, \nu} + ((\pa_t \varphi)^2 + \pa_{y_d}\varphi - t - \rho^2)a.$$

 The phase $\varphi$ satisfies the eikonal equation
$$ (\pa_t\varphi)^2 + \pa_{y_d}\varphi - t - \rho^2 - \ii \mu \sum_{|\alpha| = 1}^M (\pa_{\eta}^{\alpha} q) g_{\alpha}(\varphi) = R_M(t),$$
with $g_{\alpha}(\varphi) = \frac{1}{|\alpha|\!} \prod _{j= 1}^{n-1} (\pa_{y_j}\varphi)^{\alpha_j}$ and $R_M(t) = {\mathcal O}(t^M).$
We choose $\varphi_1 = \rho$ and one determines $\varphi_k, \: k \geq 2,$ from the equation
$$\sum_{k + j = K}(k+1)(j+1) \varphi_{k+1} \varphi_{j+1} + \pa_{y_d} \varphi_K + \epsilon_K = F(\varphi_1,...,\varphi_K)$$
with $\ep_1 = -1,\: \ep_K = 0$ for $K \geq 2$ and $F(\varphi_1,...,\varphi_K)$ given by the equality (6.6) in \cite{V1}. Next the  functions $a_{k, \nu}$ are determined form the equations
$$2 i \sum_{j=0}^{\nu} (j+1) (\nu + 1- j) \varphi_{\nu + 1 -j} a_{k, j+1} + (\nu+1)(\nu+2) a_{k-1, \nu+2} + i \pa_{y_d} a_{k, \nu} $$
$$ = \sum _{|\alpha| = 0}^M \sum_{k'= 0}^k \sum_{\nu'=0}^{\nu} b_{\alpha, k', \nu, \nu'} \pa_{y}^{\alpha} a_{k', \nu'}.$$
Therefore Lemma 6.1, 6.2, 6.3, 6.4 in \cite{V1} hold without any change since the sign before $t$ in the form of $P_0$ is not involved. Thus, as in Section 6 of \cite{V1}, for a neighbourhood $Y$ of a point in $\R^{d-1}$ obtain 
\begin{proposition} Assume $($\ref{eq:5.6}$)$ and $($\ref{eq:5.7}$)$ fulfilled  for $\eta \in {\text supp} \: \Phi_1.$ Then for all $ s \geq 0$ we have the estimates
\begin{eqnarray}
\|P_0\tilde{u}_1\|_{H^s(\R^+ \times Y)} \leq C_{s, M} h^{M\ep/2} \|f\|_{L^2(Y)}, \label{eq:5.8})\\
\|D_t \tilde{u}_1\vert_{t = 0} \|_{L^2(Y)} \leq C h^{\ep}\|f\|_{L^2(Y)}. \label{eq:5.9}
\end{eqnarray}
\end{proposition}
 To cover the region $h^{2/3} \leq  \mu \leq h^{\ep}$, it remains to study  the case when $h^{2/3} \leq \mu < h^{2/3- \ep}$ and $|\eta_d| \leq h^{1/2 - \ep}$. For these values of $\mu$ and $\eta_d$  the condition
\begin{equation} \label{eq:5.10}
\mu (\mu + |\eta_d|) \leq h^{1 + \ep}
\end{equation}
is satisfied.
We will construct a parametrix for the problem
\begin{equation} \label{eq:5.11}
\begin{cases} P_0 u = 0 \: {\rm in}\: \R^+ \times Y,\\
u = f_2 \: {\rm on}\: Y \end{cases}
\end{equation}
with $f_2 = Op_h\Bigl(\phi(\frac{\eta_d}{h^{1/2 - \ep}})\Bigr) f + \oc (h^{\infty}) f, \: f \in L^2(Y)$. For the construction we need some estimates for the Airy function $A(z) = Ai( e^{i 2 \pi/3} z)$. Here
$Ai(z)$ is the Airy function defined for $s \in \R$ by 
$$Ai(s) = \frac{1}{2\pi} \int_{-\infty}^{\infty}  e^{\ii(st + t^3 /3)} dt.$$
In the following the branch $-\pi < \arg z < \pi$ will be used and $z^{1/2} = |z|^{1/2} e^{\ii \arg z/ 2}.$ Notice also that
$$\re \sqrt{z} \geq \frac{|\im z|}{2|z|^{1/2}},\:  \im \sqrt{z}=  \frac{\im z}{2 \re \sqrt{z}}.$$
 The function $A(z)$ satisfies the
equalities 
\begin{equation} \label{eq:5.12}
(\pa_z^2 - z) A^{(k)}(z) = k A^{(k-1)} (z),\: k \in \N,
\end{equation}
where $A^{(k)}(z) = \frac{d^k A(z)}{dz^k}.$ It is well known (see \cite{Ol}, \cite{MT}) that $A(z)$ has for $|\arg z - \frac{\pi}{3} | \geq \delta> 0$ the representation
$$A(z) = \Xi(\omega z) \exp \Bigl( \frac{2}{3}\ii (-z)^{3/2}\Bigr),$$
where $\omega= e^{2 \pi \ii/3}$ and
$$\Xi(z) \sim z^{-1/4} \sum_{j = 0}^{\infty}  a_j z^{-3j/2},\: a_0 = \frac{1}{4} \pi^{-3/4},\: |z| \to \infty.$$
In the same domain in $\C$ one has also an asymptotic expansion for the derivatives of $A(z)$ by taking in the above expansion differentiation term by term (see \cite{Ol}).
Introduce the function
$$F(z) = \frac{A'(z)}{A(z)}.$$
Then for $|\arg z - \pi/3| \geq \delta > 0$ we have 
$$F(z) = z^{1/2} \sum_{k= 0}^{\infty} b_k z^{-k},\: |z| \gg 1, b_0 \neq 0.$$
For large $|z|$ and $\im z < 0$ we have the estimate $|F(z)| \leq C |z|^{1/2}$, while for bounded $|z|$ and $\im z < 0$ one obtains $|F(z)| \leq C_1$. 
Consequently, 
$$|F(z) \leq C_0 (|z| + 1)^{1/2},\: \im z < 0.$$ 
For the derivatives $F^{(k)}(z) = \frac{\pa^k F}{\pa z^k}(z)$ (see Chapter 5 in \cite{MT}) we get the following
\begin{lemma} For $\im z < 0$ and every integer $k \geq 0$ we have the estimate 
\begin{equation} \label{eq:5.13}
|F^{(k)}(z) | \leq C_k ( |z| + 1)^{1/2 - k}.
\end{equation}
\end{lemma}
 Given an integer $k \geq 0$, set $\Phi_0(z) = 1,$
$$\Phi_k(z): = A(z) \pa_z^k(A(z)^{-1}) = \pa_z \Phi_{k-1}(z) - F(z) \Phi_{k-1}(z),\: k \geq 1.$$
 Taking the derivatives in the above equality and using (\ref{eq:5.13}), by induction in $k$ one obtains 
\begin{lemma} For $\im z < 0$ and all integers $k \geq 1, l \geq 0$, we have the bound
\begin{equation} \label{eq:5.14}
|\pa_z^{l} \Phi_k(z)| \leq C_{k, l} \Bigl( |z| + 1\Bigr)^{\frac{k }{2}- l}.
\end{equation}
\end{lemma}
For $t \geq 0$ and $\im z < 0$, set
$$\Psi_k(t, z): = \frac{A^{(k)}(-t + z)}{A(z)},\: \Psi_k^{(l)}(t, z): = \pa_z \Psi_k(t, z).$$
The next Lemma is an analogue of Lemma 3.3 in \cite{V1}.
\begin{lemma} For $\im z < 0$ and all integers $k \geq 0, l \geq 0$, we have the estimate
\begin{equation} \label{eq:5.15}
|\Psi_k^{(l)}(0, z)| \leq C_k |\im z|^{-l} \Bigl( |z|^{1/2} + 1\Bigr)^{k}.
\end{equation}
For $0 < t \leq |z| ,\: \im z < 0$ and all integers $k \geq 0, l \geq 0$, we have 
\begin{equation} \label{eq:5.16}
|\Psi_k^{(l)}(t, z)| \leq C_{k,l} |\im z|^{-l} (|z|^{1/2} + |\im z|^{-1})\Bigl( |z|^{1/2} + 1 \Bigr)^{k},
\end{equation}
while for $|t| \geq |z|$ one obtains
\begin{equation} \label{eq:5.17}
|\Psi_k^{(l)}(t, z)| \leq C_{k,l} |\im z|^{-l}(|z|^{1/2} + |\im z|^{-1}) \Bigl( t^{1/2} + |\im z|^{-1} \Bigr)^{k} e^{-t^{1/2} |\im z|/4}.
\end{equation}
\end{lemma}

{\it Proof.} Since $\Psi(t, z)$ is analytic for $\im z < 0$, it is sufficient to establish the above estimates for $l = 0$ and to apply Cauchy formula for the derivatives (see Section 3 in \cite{V1}). 
Taking into account (\ref{eq:5.12}), (\ref{eq:5.13}), by induction in $k$ one deduces
$$|A^{(k)}(z)| \leq C_k \Bigl(|z|^{1/2} + 1\Bigr)^k |A(z)|$$
hence
\begin{equation} \label{eq:5.18}
|\Psi_k(t, z)| \leq C_k \Bigl( t^{1/2} + |z|^{1/2} + 1\Bigr)^k |\Psi_0(t, z)|.
\end{equation}
Thus it is sufficient to estimate $|\Psi_0(t, z)|.$
 The representation of $A(z)$ with phase $\exp\Bigl(( \frac{2}{3}\ii(-z)^{3/2}\Bigr)$ mentioned above holds for $\im z < 0.$ Hence
$$\Bigl|\frac{A(-t + z)}{A(z)}\Bigr| \leq \Bigl|\frac{\Xi(\omega(-t + z))}{\Xi(\omega z)}\Bigr|\exp\Bigl(- \im \frac{2}{3} \Bigl((t-z)^{3/2} - (-z)^{3/2}\Bigr)\Bigr)$$
$$ = \Bigl|\frac{\Xi(\omega(-t + z))}{\Xi(\omega z)}\Bigr|e^{-\varphi}.$$
It clear that $|\Xi(\omega(-t + z))|\leq c_0$. For $|z| \leq C,\: C \gg 1$ we have 
$$\big|\Bigl(\Xi(\omega z)\Bigr)^{-1}\big|\leq C_1 \leq C_2 |\im z|^{-1},$$
 while for  $|z| \geq C$ we have 
$$\big|\Bigl(\Xi(\omega z)\Bigr)^{-1}\big| \leq C_3 |z|^{1/4} \leq C_3 |z|^{1/2}.$$
Thus $|\frac{\Xi(\omega(-t + z))}{\Xi(\omega z)}| \leq C(|z|^{1/2} + |\im z|^{-1}).$

Next,  we get
$$\varphi = \frac{2}{3} \im (t - z)^{3/2} - \frac{2}{3}\im (-z)^{3/2} = \int_0^t \im (\tau - z)^{1/2} d\tau.$$
$$ = -\int_0^t \frac{\im z}{2\re((\tau-z)^{1/2})}d\tau \geq \frac{t |\im z|}{2( t^{1/2} + |z|^{1/2})}$$
and this shows that for $t > 0$ we have $\varphi > 0.$
For $|t| \leq |z|$ the estimate (\ref{eq:5.18}) implies (\ref{eq:5.16}).
For $|t| \geq |z|$ we have 
$$\frac{t |\im z|}{2( t^{1/2} + |z|^{1/2})}\geq \frac{t^{1/2} |\im z|}{4}$$
and 
\begin{equation} \label{eq:5.19}
t^{k/2} e^{-t^{1/2} |\im z|/4} \leq C_k |\im z|^{-k}.
\end{equation}
If $|\im z| \leq 1$ we have $1 \leq |\im z|^{-1}$, while if $|\im z| > 1$, we get $t \geq |z| > 1$. Hence from (\ref{eq:5.19}) and (\ref{eq:5.18}) we deduce (\ref{eq:5.17}).
\hfill\qed\\

\def\ac{{\mathcal A}}
For $h^{2/3} \leq \mu \leq h^{2/3 - \ep}$ we will construct a parametrix for (\ref{eq:5.11}) repeating without any change the construction in Section 5 of \cite{V1}. The parametrix has the form $\tilde{u}_2 = \phi(t/h^{\ep}) Op_h(\ac(t)) g$, where $g \in L^2(Y)$ can be determined as in Section 5, \cite{V1}. Here
$$\ac(t) =\sum_{k= 0}^M a_k(y, \eta; h, \mu) \psi_k(t, y; h, \mu),$$
$$ \psi_k = h^{k/3} \Psi_k\Bigl(-t h^{-2/3}, (\eta_d - \ii\mu q(y, \eta)) h^{-2/3}\Bigr),$$ 
$M \gg 1$ is an arbitrary integer, $a_0 = \phi(\frac{\eta_{d}}{h^{1/2- \ep}})$. Next $a_k, \: k \geq 1$, are independent on $t$ and they can be determined from the equality
$$(k + 1) a_{k+ 1} = - \ii \pa_{y_d} a_k + \mu h^{-2/3} \pa_{y_d}q F(\eta_d - \ii \mu q(y,\eta) h^{-2/3}) a_k - \mu h^{-1} \pa_{y_d} a_{k-1}$$
$$+ \sum_{l = 0}^{k} \sum_{|\alpha| = 0}^{k} \Bigl( b_{k,l, \alpha}^{(1)} + b_{k, l, \alpha}^{(2)}\Bigr) \pa_y ^{\alpha} a_k.$$ 

 We have
$$P_0 Op_h(\ac(t)) = Op_h\Bigl((D_t^2 - t + \eta_{d} - \ii \mu q(y, \eta) - \ii h \pa_{y_d}) \ac(t)\Bigr)$$
$$- \ii \mu q(y, D_y) Op_h(\ac(t)) + \ii \mu Op_h(q \ac(t)) + h \tilde{q}(y, D_y) Op_h (\ac(t)).$$
On the other hand, (\ref{eq:5.12}) implies the equality
$$(D_t^2 - t + \eta_{d} - \ii \mu q(y, \eta)) \Psi_k \Bigl(-t h^{-2/3}, (\eta_d - \ii\mu q(y, \eta)) h^{-2/3}\Bigr)$$
$$ = - kh^{-2/3} \Psi_{k-1} \Bigl(-t h^{-2/3}, (\eta_d - \ii\mu q(y, \eta)) h^{-2/3}\Bigr)$$
and 
$$(D_t^2 - t + \eta_{d} - \ii \mu q(y, \eta))\ac(t) = - h \sum_{k= 0}^{M-1} (k+1) a_{k+1} \psi_k.$$
Next the construction of the parametrix goes without any changes as in Section 5 in \cite{V1} applying Lemmas 5.5, 5.6 and 5.7 instead of Lemmas 3.1, 3.2 and 3.3 in \cite{V1}. Thus as an analogue of Theorem 5.7 in \cite{V1} we get the following
\begin{proposition} For all $s \geq 0$, we have the bounds
\begin{eqnarray}
\|P_0 \tilde{u}_2\|_{H^s(\R^+ \times Y)} \leq C_{s, M} h^{M\ep/2} \|f\|_{L^2(Y)},\\
\|\tilde{u}_2\vert_{t = 0} - Op_h\Bigl(\phi(\eta_d \mu /h^{1+ \ep})\Bigr) f \|_{L^2(Y)} \leq \oc (h^{\infty}) \|f\|_{L^2(Y)},\\
\|D_t \tilde{u}_2 \vert_{t = 0}\|_{L^2(Y)} \leq C h^{\ep} \|f\|_{L^2(Y)}.
\end{eqnarray}
\end{proposition}
Combining Proposition 5.4 and Proposition 5.8, we obtain, as in \cite{V1}, Theorem 5.2.\\

After this preparation we pass to the analysis of an eigenvalues-free region when 
$$\re z = 1,\:h^{2/3} \leq \im z \leq h^{\ep},\: 0 < \ep \ll 1.$$
Let $\rho = \sqrt{ 1 - r_0 + \ii \im z}$. As in the previous section, we examine the equation
$$\nc_{ext}(z, h)(f) - \sqrt{z} \gamma f = 0.$$
Consider the partition of the unity $\chi_{\ep/2}^+ + \chi_{\ep/2}^0 + \chi_{\ep/2}^- = 1$ on $T^*(\Gamma)$ introduced in the beginning of this section. Applying Theorem 5.2, we have
$$\|\nc_{ext}(z, h) ( 1- \chi_{\ep/2}^0) f - \sqrt{z} \gamma f \|_{\lg} \leq C h^{\ep/2} \|f\|_{\lg}.$$
Taking into account Theorem 5.1 for the operators $N(z, h) \chi_{\ep/2}^{\pm}$, one deduces
\begin{equation} \label{eq:5.21}
\|Op_h\Bigl(\rho (\chi_{\ep/2}^+ + \chi_{\ep/2}^-) - \sqrt{z}\gamma  \Bigr)f\|_{\lg} \leq C_1 h^{\ep/2} \|f\|_{\lg}.
\end{equation}
We write
$$g_1: = \rho (\chi_{\ep/2}^+ + \chi_{\ep/2}^-) - \sqrt{z}\gamma = \frac{\rho^2 \Bigl[(\chi_{\ep/2}^+)^2 + (\chi_{\ep/2}^-)^2\Bigr] - z \gamma^2}{
\rho (\chi_{\ep/2}^+ + \chi_{\ep/2}^-) + \sqrt{z}\gamma}.$$
Clearly,
$$\re \Bigl(\rho^2 \Bigl[(\chi_{\ep/2}^+)^2 + (\chi_{\ep/2}^-)^2\Bigr] - z \gamma^2\Bigr) = (1 - r_0)\Bigl[(\chi_{\ep/2}^+)^2 + (\chi_{\ep/2}^-)^2\Bigr]- \gamma^2\leq - \eta_0 < 0$$
since $1 - r_0 \leq 1,$ supp $\chi_{\ep/2}^+ \cap \:{\rm supp}\: \chi_{\ep/2}^{-} = \emptyset$ and $1 - \gamma^2 \leq - \eta_0.$
Thus for bounded $|\xi'|$ we have $|g_1| \geq \eta_2 > 0$, while for $|\xi'| \gg 1$ we get $|g_1| \sim |\xi'|$. To estimate $g_1^{-1}$, it is necessary to estimate only $\rho (\chi_{\ep/2}^+ + \chi_{\ep/2}^-) + \sqrt{z}\gamma$ and one deduces
$$|\pa_{x'}^{\alpha} \pa_{\xi'}^{\beta} (g_1^{-1})|\leq C_{\alpha, \beta} h^{-\frac{\ep}{2} ( 1/2 +|\alpha| + \beta|)} ( 1 + |\xi'|)^{1 - |\beta|}.$$
The same estimates holds for $g_1$, hence $g_1 \in S_{\ep/2}^{\ep/4, 1},\: g_1^{-1} \in S_{\ep/2}^{\ep/4, -1}$  and 
$$\|(Op_h(g_1^{-1}) Op_h(g_1) - I) f\|_{\lg} \leq C h^{1 - \ep} \|f\|_{\lg}. $$
Combining this with (\ref{eq:5.21}), for small $h$ we conclude as in Section 4, that $f = 0.$\\

It remains to study the case 
$$z \in {\mathcal D} = \{ z \in \C: z = 1 + \ii \im z, \: 0 < \im z \leq h^{2/3}\}.$$
 The  Dirichlet problem for $-h^2 \Delta - z$ with $z = 1 + \ii \im h^{2/3} w,\: |w| \leq C_0,$ has been investigated by  Sj\"ostrand in Chapters 9 and 10 in \cite{Sj} (see also \cite{CPV}). For $0 \leq w \leq 1$ this covers the region ${\mathcal D}$. In \cite{Sj} the exterior Dirichlet-to-Neumann map ${\mathcal N}_{ext}(z, h)$ is defined with respect to the outgoing solution \footnote{the outgoing solutions in the sense of Lax-Phillips \cite{LP1} are different from the outgoing ones in \cite{Sj}. See Section 6 for more details.} of the problem (\ref{eq:1.5}). Notice that for $\im z > 0$ the outgoing solutions are in $H_h^2(\Omega)$, so the exterior Dirichlet-to-Neumann map in \cite{Sj} coincides with that defined in Section 2. We need to recall some results in Chapter 10 of \cite{Sj}.  The operator $\nc_{ext}(z, h)$ is a $h-$pseudo-differential operator with symbol $n_{ext}(x', \xi', h)$. Introduce the glancing set
$${\mathcal G} = \{(x', \xi') \in T^*(\Gamma):\: r_0(x', \xi') = 1\}.$$
We have $\gamma(x) \geq 1 + \eta_3 > 1,\: \forall x \in \Gamma.$ Choose a small number $\delta_0,\:0 < \delta_0 < \eta_3/2$. Then for  $|r_0(x', \xi') - 1| \geq \delta_0$ the symbol $n_{ext}$ satisfies the estimates
\begin{equation} \label{eq:5.24}
|\pa_{x'}^{\alpha} \pa_{\xi'}^{\beta} n_{ext}(x', \xi', h)| \leq C_{\alpha, \beta}  \langle \xi'\rangle^{1 -|\beta|},\: \forall \alpha, \forall \beta,
\end{equation}
while for $|r_0(x', \xi') -1| \leq 2\delta_0$ we have the estimates 
\begin{equation} \label{eq:5.25}
|\pa_{x'}^{\alpha} \pa_{\xi'}^{\beta} n_{ext}(x', \xi', h)| \leq C_{\alpha, \beta} (h^{2/3} + |r_0 - 1|)^{\frac{1}{2} - \beta_d},\: \forall\alpha, \forall \beta
\end{equation}
if $r_0(x', \xi') - 1$ is transformed into $\xi_{d}$
 by a tangential Fourier integral operator as it was mentioned in the beginning of this section. From the estimates near ${\mathcal G}$ it follows that for small $0 < h \leq h_0(\delta_0)$ we have a bound 
$$ \Bigl\|\nc_{ext}(z, h) \phi\big(\frac{1 - r_0(x', \xi')}{\delta_0}\big)\Bigr\|_{\lg \to \lg} \leq C (h^{1/3} + \delta_0^{1/2})$$ 
with a constant $C > 0$ independent on $h$ and $\delta_0.$ Let $f \neq 0$ be the trace of $u\vert_{\Gamma}$,  where $(u, v)$ is an eigenfunction of $G$. Consider the equality
\begin{eqnarray}
 -\re \Bigl\langle \nc_{ext}(z, h) \Bigl[ 1 - \phi\big(\frac{1 - r_0(x', \xi')}{\delta_0}\big)\Bigr] f, f \Bigr\rangle _{\lg} +  \re \langle \sqrt{z}\gamma f, f \rangle_{\lg} \nonumber \\
 = \re \langle \nc_{ext}(z, h)  \phi\big(\frac{1 - r_0(x', \xi')}{\delta_0}\big) f, f \rangle_{\lg}. \label{eq:5.26}
\end{eqnarray}
The above estimate shows that the right hand side in (\ref{eq:5.26}) is bounded by  $C_1 (h^{1/3} + \delta_0^{1/2})\|f\|^2_{\lg}.$\\

 Introduce two functions
$\psi_{\pm}(\sigma) \in C^{\infty}(\R:[0, 1])$ such that $\psi_{+}(\sigma) = 0$ for $\sigma \leq 1/2,\: \psi_{+}(\sigma) = 1$ for $\sigma \geq 1, \: \psi_{-}(\sigma) = \psi_{+}(-\sigma).$ We write
$$\nc_{ext}(z, h) \Bigl[ 1 - \phi\big(\frac{1 - r_0(x', \xi')}{\delta_0}\big)\Bigr] = \nc_{ext}(z, h) \chi_{+} + \nc_{ext}(z, h) \chi_{-},$$
where 
$$\chi_{\pm}(x', \xi') = \Bigl[ 1 - \phi\big(\frac{1 - r_0(x', \xi')}{\delta_0}\big)\Bigr]\psi_{\pm}\Bigl( \frac{1 - r_0(x', \xi')}{\delta_0}\Bigr)$$
have support in 
$\{ (x', \xi'): 1 - r_0(x', \xi') \geq \delta_0/2\}$ and 
 $\{ (x', \xi'): 1 - r_0(x', \xi') \leq - \delta_0/2\}$, respectively. The principal symbols $n_{\pm}$ of $\nc_{ext}(z, h) \chi_{\pm}$ have the form
$$n_{\pm}= \Bigl(\sqrt{1 - r_0 + \ii \im z} \Bigr)\chi_{\pm}$$
and 
$$ \re \langle \nc_{ext}(z, h) \chi_{\pm} f, f \rangle_{\lg} = \langle Op_h(\re (n_{\pm})) f, f \rangle_{\lg} + {\mathcal O}(h)\|f\|_{\lg}^2.$$
On the other hand, 
$$|\re n_{+}| = \chi_{+} |\re \sqrt{1 - r_0 + \ii \im z}| \leq (1 + h^{2/3})^{1/2}.$$
In the same way for the principal symbol $n_{-}$ of $\nc_{ext}(z, h) \chi_{-}$ we get
$$|\re n_{-}| = \chi_{-} |\re \sqrt{1 - r_0 + \ii \im z}| \leq  y^{1/2} \sin \frac{\psi}{2},$$
where $1 - r_0 + \ii \im z = y e^{\ii(\pi - \psi)},\: y > 0, \: 0 < \psi \ll 1.$ Next
$$ y^{1/2} \sin \frac{\psi}{2} = y^{1/2}\sqrt{ \frac{1 - \cos \psi}{2}} = \frac{1}{\sqrt{2}} \sqrt{ y -(r_0- 1)}.$$
On the support of $\chi_{-}$ we have $0 < r_0 - 1 < y \leq (r_0 - 1) + h^{2/3}$, and this implies $ y^{1/2} \sin \frac{\psi}{2} \leq \frac{1}{\sqrt{2}} h^{1/3}.$ Combining the above estimates, we conclude that
\begin{equation} \label{eq:5.27}
 -\re \Bigl\langle \nc_{ext}(z, h) \Bigl[ 1 - \phi\big(\frac{1 - r_0(x', \xi')}{\delta_0}\big)\Bigr] f, f \Bigr\rangle_{\lg} \geq -(1 + C_1h^{1/3})\|f\|_{\lg}^2.
\end{equation}
Let $ \sqrt{z} = v + \ii w.$ Then $v^2 = 1 + w^2 \geq 1$ yields $\re \sqrt{z} = v \geq 1.$ Consequently,
$$\re \langle \sqrt{z} \gamma f, f \rangle_{\lg} \geq (1 + \eta_3) \|f\|^2_{\lg}.$$
From this estimate and (\ref{eq:5.27}) one deduces that the left hand side of (\ref{eq:5.26}) is greater than $(\eta_3 - C_1 h^{1/3})\|f\|^2_{\lg}$. For small $h$ and small $\delta_0$ (depending on $\eta_3$) we obtain a contradiction with the estimate of the right hand side of (\ref{eq:5.26}). Finally, if $\re z = 1, 0 <  \im z \leq h^{2/3}$ with $0 < h \leq h_0(\eta_3)$ there are no eigenvalues $\lambda = \frac{\ii \sqrt{z}}{h}$ of $G$.  Combining this with the result of Section 4, completes the proof of Theorem 1.3.\\

\section{Trace formula}

Before going to the proof of a trace formula for the counting function of the eigenvalues of $G$, we need to examine the properties of the Dirichlet-to-Neumann map $\nc(\lambda)$ defined below.
This map can be used to prove the discreteness of the spectrum of $G$ in $\{z \in \C: \re z < 0\}$. This result for $d$ odd was established in \cite{LP}  and the proof there exploits the fact that
the scattering operator $S(z)$ is invertible for $z = 0.$ For even dimensions $d$ this property of $S(z)$ is not true. 
We present a proof of the discreteness of the spectrum of $G$ based on the invertibility of an operator involving $\nc(\lambda)^{-1}$ and it seems that for $d$ even this result is new.

\begin{proposition} Let $\gamma(x) \neq 1$ for all $x \in \Gamma$. Then for $d \geq 2$ the spectrum of the generator $G$ in $\{z \in \C:\: \re z < 0\}$ is formed by isolated eigenvalues with finite multiplicities.
\end{proposition}

{\it Proof.}  Consider for $\re \lambda < 0$ the  map 
$$\nc(\lambda): H^s(\Gamma) \ni f \longrightarrow \pa_{\nu} u\vert_{\Gamma} \in H^{s-1}(\Gamma),$$
where $u$ is the solution of the problem
\begin{equation} \label{eq:6.1}
\begin{cases} (\Delta - \lambda^2) u = 0 \: {\rm in}\: \Omega,\: u \in H^2(\Omega),\\
u = f \:{\rm on}\: \Gamma. \end{cases}
\end{equation}

The condition $u \in H^2(\Omega)$ implies that $u$ is $\ii\lambda$- outgoing which means that there exists $R > \rho_0$ and a function $g \in L^2_{comp}(\R^d)$ such that
$$u(x) = (-\Delta_0 + \lambda^2)^{-1} g,\: |x| \geq R,$$
where $R_0(\lambda) =(-\Delta_0 +\lambda^2)^{-1}$ is the outgoing resolvent of the free Laplacian in $\R^d$ which is analytic for $\re \lambda > 0.$ Recall that $R_0(\lambda)$ has kernel

 \be \label{eq:6.2}  R_0(\lambda, x- y) = -\frac{\ii}{4} \Bigl(\frac{-\ii \lambda}{2 \pi |x-y|}\Bigr)^{(n-2)/2}\Bigl( H_{\frac{n-2}{2}}^{(1)}(u)\Bigr) \Bigl\vert_{u = -\ii \lambda |x-y|},
\ee
where $H_k^{(1)}(x)$ is the Hankel function of first kind and
we have the asymptotic (see for example, \cite{Ol})

\be \label{eq:6.3}
H^{(1)}_{\mu}(z) = \Bigl(\frac{2}{\pi r}\Bigr)^{1/2} e^{\ii (z - \frac{\mu \pi}{2} - \frac{\pi}{4})} + {\mathcal O}(r^{-3/2}), \: - 2\pi < {\rm arg} z < \pi,|z| =r  \to +\infty
\ee

Below we present some well known facts for the sake of completeness. The solution of the Dirichlet problem
(\ref{eq:6.1}) with $f \in H^{3/2}(\Gamma)$ has the representation
$$u = e(f) + (-\Delta_D +\lambda^2)^{-1}(\Delta -\lambda^2) (e(f)),$$ 
where $e(f): H^{3/2}(\Gamma) \ni f \to e(f) \in H^{2}_{comp}(\Omega)$ is an extension operator  and $R_D(\lambda) =(-\Delta_D +\lambda^2)^{-1}$ is the outgoing resolvent of the Dirichlet Laplacian $\Delta_D$ in $\Omega$ which is analytic for $\re \lambda < 0$.\footnote{Notice that the definition of outgoing solutions in \cite{LP1} is different from that given above  and the outgoing solutions in our paper correspond to incoming ones in \cite{LP1}. To avoid misunderstanding the precise form of $R_0(\lambda, x-y)$ is given in (\ref{eq:6.2}).}   Therefore
$$\nc(\lambda) f= \pa_{\nu}(e(f)) + \pa_{\nu} \Bigl[(-\Delta_D +\lambda^2)^{-1}(\Delta -\lambda^2)(e(f))\Bigr]$$
implies that $\nc(\lambda)$ is analytic for $\Re \lambda < 0$.  The solution of (\ref{eq:6.1}) for $\Re \lambda < 0$ can be written also as follows
$$u(x; \lambda) = \int_{\pa \Omega} \Bigl[ R_0( \lambda , x- y)(N(\lambda) f)(y) - 
\frac{\pa R_0( \lambda, x - y)}{\pa \nu_y} f(y) dy\Bigr]. $$
Taking the trace on $\Gamma$,  this implies
$$ C_{00}(\lambda) + C_{01}(\lambda) N(\lambda) =  Id,$$
where $C_{00}(\lambda)$ and $C_{01}(\lambda)$ are the Calderon operators (see for example, \cite{Me}) which are analytic operator-valued functions for $\lambda \in\C$ for $d$ odd and on the logarithmic covering of $\C$ for $d$ even. Melrose proved (\cite{Me}, Section 3) that there exists an entire family $P_D(\lambda)$ of pseudo-differential operators of order -1 on $\Gamma$ so that
$$(-\Delta_{\Gamma} + 1)^{1/2} C_{01}(\lambda) = Id + P_D(\lambda),$$
$\Delta_{\Gamma}$ being the Laplace Beltrami operator on $\Gamma$. For $\re \lambda < 0$ this implies
$$N(\lambda) = (Id + P_D(\lambda))^{-1} (-\Delta_{\Gamma} + 1)^{1/2} (Id - C_{00}(\lambda)).$$
On the other hand, it is well known that the Neumann problem 
\begin{equation} \label{eq:6.4}
\begin{cases} (\Delta - \lambda^2) u = 0 \: {\rm in}\: \Omega,\: u \in H^2(\Omega),\\
\pa_{\nu}u = 0 \:{\rm on}\: \Gamma. \end{cases}
\end{equation}
 has no non trivial $(\ii\lambda)$-outgoing solutions for $\re \lambda < 0.$ This implies that for $\re \lambda < 0$ the operator  $C_{00}(\lambda)$ has not 1 as an eigenvalues  and since $C_{00}(\lambda)$ is compact by the analytic Fredholm theorem we deduce that
 $N(\lambda)^{-1}$ is analytic for $\re \lambda < 0.$\\

Going back to the problem (\ref{eq:1.2}), we write the boundary condition  as follows
$$\nc(\lambda) \Bigl( I - \lambda \nc^{-1}(\lambda) \gamma \Bigr) f_1 = 0,\:\re \lambda < 0,\: x \in \Gamma.$$
The operator $\nc(\lambda)^{-1}: L^2(\Gamma) \longrightarrow H^{1}(\Gamma)$ is compact and by Theorems 1.2 and Theorem 1.3 there are points $\lambda_0,\: \re \lambda_0 < 0,$ for which $(I - \lambda_0 \nc^{-1}(\lambda_0)\gamma)$ is invertible. Applying the analytic Fredholm theorem, one concludes that the spectrum of $G$ in the open half-plane $\re \lambda < 0$ is formed by isolated eigenvalues with finite multiplicities.\\

\begin{remark} The assumption  $\gamma(x) \neq 1,\: \forall x \in \Gamma,$ was used only to apply Theorems $1.2$ and $1.3$. For odd dimensions $d$ we can relax this assumption. Indeed, for $d$ odd we have no resonances in a small neighbourhood of 0 for the Dirichlet and Neumann problems, so we may apply the above argument in a open domain including a small neighbourhood of 0. For $d$ even this property does not hold. \footnote{In \cite{M1} one obtains eigenvalues-free regions in the case $\gamma \geq 1$, but in this paper  one applies the result of \cite{LP} for $d$ odd.} 
\end{remark}

Now we pass to a trace formula involving the operator 
$$ C(\lambda) := \nc(\lambda) - \lambda \gamma = \nc(\lambda) \Bigl( I - \lambda \nc^{-1}(\lambda) \gamma \Bigr),$$
which by the analysis above is an analytic operator-valued function in $\{ z \in \C: \: \re \: z < 0\}$, while $C(\lambda)^{-1}$ is meromorphic in the same domain.
 Our purpose is to find  a formula for the trace
\begin{equation} \label{eq:6.5}
  {\rm tr}\: \frac{1}{2 \pi i} \int_{\delta} (\lambda - G)^{-1} d\lambda,
\end{equation}
 where $\omega \subset \{\re\: z < 0\}$ is a domain with boundary  the a positively oriented curve $\delta$ and $(G-\lambda)^{-1}$ is analytic on $
\delta$. Since $(G- \lambda)^{-1}$ is meromorphic in $\omega$,  if $\lambda_0$ is a pole of $(G- \lambda)^{-1},$ the (algebraic) multiplicity of an eigenvalue $\lambda_0$ of $G$ is given by
$${\rm mult}\: (\lambda_0) ={\rm rank} \frac{1}{2 \pi i} \int_{|\lambda - \lambda_0| = \ep_0} (\lambda - G)^{-1} d\lambda,$$
with $\ep_0 > 0$ small enough and $\{\lambda \in \C:\:|\lambda - \lambda_0| = \ep_0\}$ positively oriented. Therefore, (\ref{eq:6.5}) is just equal to the number of the eigenvalues of $G$ in $\omega$ counted with their multiplicities.
Let $(u, w) =  (G-\lambda)^{-1} (f, g)$. Then we have $w = \lambda u + f$ and 
$$u = -R_D(\lambda) (g + \lambda f) + K(\lambda) q.$$
Here $R_D(\lambda)= (-\Delta_D + \lambda^2)^{-1}$ is the outgoing resolvent introduced in the proof of Proposition 6.1 and $K(\lambda)$ satisfies
$$\begin{cases}  (\Delta - \lambda^2) K(\lambda) = 0 \:\: {\rm in}\: \Omega,\\
K(\lambda) = \:Id \:\:{\rm on}\: \Gamma.
\end{cases}$$
The boundary condition on $\Gamma$ implies
$$\pa_{\nu} \Bigl[-R_D(\lambda) ( g + \lambda f) + K(\lambda) q\Bigr]- \gamma \lambda \Bigl[- R_D(\lambda) ( g + \lambda f) +  q\Bigr] - \gamma f = 0,\: x \in \Gamma$$
and the term $\gamma \lambda [R_D(\lambda) ( g + \lambda f) $ vanishes. Next $\nc(\lambda) = \pa_{\nu} K(\lambda)\vert_{\Gamma}$ is the Dirichlet-to-Neumann map, and assuming  $C^{-1}(\lambda)$ invertible, one gets
 $$q = C^{-1}(\lambda) \Bigl( [ \pa_{\nu} R_D(\lambda)(g + \lambda f) ] + \gamma f \Bigr).$$
Therefore
$$u = \Bigl[-\lambda R_D(\lambda) + \lambda K(\lambda) C^{-1}(\lambda) \pa_{\nu} R_{D}(\lambda) + C^{-1}(\lambda)\gamma\Bigr] f  + Xg$$
$$ w = Y f + \Bigl[\lambda R_{D}(\lambda) + \lambda K(\lambda) C^{-1}(\lambda) \pa_{\nu} R_D(\lambda) \Bigr]g,$$
where the form of the operators $X$ and $Y$ is not important for the calculus of the trace.
Thus we have the equality
$${\rm tr}\:\int_{\delta} (\lambda - G)^{-1} d\lambda = -{\rm tr}\: \int_{\delta} \Bigl(2 \lambda K(\lambda) C^{-1}(\lambda) \pa_{\nu} R_D(\lambda) + C^{-1}(\lambda) \gamma\Bigr) d\lambda.$$
The operator $C^{-1}(\lambda)$ is meromorphic with finite rang singularities near every pole. For the first term in the integral on the right hand side we apply Lemma 2.2 in \cite{SjV} combined with the fact
$$\frac{\pa \nc}{\pa \lambda}(\lambda)= \pa_{\nu} \frac{\pa K}{\pa \lambda}(\lambda)\vert_{\Gamma}  = - 2 \lambda \pa_{\nu}R_{D} (\lambda) K(\lambda).$$
Finally, we obtain the following
\begin{proposition} Assume $\gamma(x) \neq 1, \forall x \in \Gamma.$ Let $\delta \subset \{ z \in \C: \: \re\: \lambda < 0\}$ be a closed positively oriented curve and let $\omega$ be the domain bounded by $\delta$. Assume that  $C^{-1}(\lambda)$ is meromorphic in $\omega$ without poles on $\delta$ . Then
\begin{equation} \label{eq:6.6}
 {\rm tr}\: \frac{1}{2 \pi i} \int_{\delta} (\lambda - G)^{-1} d\lambda = {\rm tr} \frac{1}{2 \pi i} \int_{\delta} C^{-1}(\lambda) \frac{\pa C}{\pa \lambda}(\lambda) d \lambda.
\end{equation}
\end{proposition}

In the case $(B)$ it is interesting to apply  Proposition 6.3  to obtain a Weyl formula for the eigenvalues of $G$ lying in the domain ${\mathcal R}_N$ following the approach in \cite{SjV} and \cite{PV}.

\section{Appendix}

\def\Re{{\rm Re}\:}
\def\Im{{\rm Im}\:}
\def\ep{\epsilon}
\def\be{\begin{equation}}
\def\ee{\end{equation}}
\def\cn{{\mathcal N}}
\def\sn{{\mathbb  S}^{n-1}}
\def\Ker {{\rm Ker}\:}
\def\la{\langle}
\def\ra{\rangle}

I this Appendix we assume that $\gamma \geq 0$ is a constant and $d$ odd. We examine the existence of the eigenvalues of $G$ for the ball $B_3 = \{x \in \R^3: |x| \leq 1\}$. 
Consider the Dirichlet problem for the Helmholtz equation in the exterior of $B$. 
\be \label{eq:7.1} \begin{cases}  (\Delta - \lambda^2) u = 0\: {\rm in}\: |x| > 1,\: u \in H^2(|x| \geq 1),\cr
u\vert_{|x| = 1} = f \in L^2(\S^2). \end{cases}
\ee
Setting $\lambda = i \mu$, $\Im \mu > 0,$ it is well known that the outgoing solution of (\ref{eq:7.1}) in polar coordinates $(r, \omega), r \in \R^+,\:\omega \in \S^2$ is given by a series
$$u(x, \mu) = \sum_{n = 0}^{\infty}\sum_{m=-n}^{n} a_{n, m} \frac{h^{(1)}_{n}(\mu r)}{ h^{(1)}_{n}(\mu) }Y_{n, m}(\omega), \: |x| = r.$$
Here $Y_{n, m}(\omega)$ are the spherical functions which are eigenfunctions of the Laplace-Beltrami operator $-\Delta_{\S^2}$ with eigenvalues $n(n+1)$ and 
$$h^{(1)}_{n}(r) = \frac{H^{(1)}_{n + 1/2}(r)}{r^{1/2}}$$ are the spherical (modified) Hankel functions of first kind. The boundary condition in (\ref{eq:7.1}) is satisfied choosing $a_{n,m}$ so that
\be \label{eq:7.2}
f(\omega) = \sum_{n=0}^{\infty}\sum_{m = - n}^{n} a_{n,m} Y_{n, m} (\omega).
\ee
Now consider the boundary problem
\be \label{eq:7.3} \begin{cases}  (\Delta - \lambda^2) u = 0\: {\rm in}\: |x| > 1, u \in H^2(|x| \geq 1),\cr
\partial_r u - \lambda \gamma u = 0\: {\rm on}\:\S^{2}. \end{cases}
\ee

We will prove the following
\begin{proposition}
For $\gamma = 1$ and $\Re \lambda < 0$ there are no non trivial solutions of $(7.3)$. For $0 < \gamma < 1$ the eigenvalues of $G$ 
lie in the region
\be \label{eq:7.4} 
\Bigl\{ \lambda \in \C: \pi/4 < |\pi - {\rm arg}\:\lambda | < \pi/2,\: |\lambda| > \frac{ \cos ( {\rm arg}\: \lambda)}{(1- \gamma)\cos (2{\rm arg}\:\lambda)}\Bigr\}.
\ee

\end{proposition}

{\it Proof.} Introduce the Dirichlet-to-Neumann map
$N(\lambda)f = \partial_r u\vert_{|x| = 1},$ where $u$ is the solution of (\ref{eq:7.1}). Assume that $(u, v)$ is an eigenfunction of $G$. Then $u$ satisfies (\ref{eq:7.3}). Setting $u\vert_{|x| = 1} = f$, $\lambda = \ii \mu,\: \Im \mu > 0$, the boundary condition  implies

\be\label{eq:7.5}
N(\ii \mu)f - \ii \mu \gamma f = 0\: {\rm on} \: |x| = 1,
\ee
and we deduce 
$$\sum_{n= 0}^{\infty} \sum_{m = -n}^{n} C(n ; \mu, \gamma) a_{n, m} Y_{n. m}(\omega) = 0$$
with
\be \label{eq:7.6}
 C( n; \mu, \gamma): = \partial_r \Bigl[\frac{h^{(1)}_{n}(\mu r)}{ h^{(1)}_{n}(\mu)}\Bigr]\vert_{r = 1} - i \mu \gamma
 ,\: n \in \N.
\ee
It is well known (see \cite{Ol}) that $h^{(1)}_n(x)$ have the form
$$h^{(1)}_n(x) = (-\ii)^{n+1} \frac{e^{ix}}{x} \sum_{m=0}^n \frac{\ii^m}{m!(2x)^m} \frac{(n+m)!}{(n-m)!} = (-\ii)^{n+1} \frac{e^{ix}}{x} R_n(x).$$

The problem is reduced to show that  $C(n; \mu, \gamma) \neq 0$ for all $n \in \N$, that is
\be \label{eq:7.7}
\partial_r \Bigl( h^{(1)}_{n} (\mu r)\Bigr)\vert_{r = 1} - \ii \mu \gamma h^{(1)}_{n}(\mu) \neq 0, \: \forall n \in \N.
\ee
In fact this implies that all coefficients $a_{n,m}$ vanish , so $f = 0.$
Taking the derivative with respect to $r$, one obtains 
\be \label{eq:7.8}
C(n; \mu, \gamma) = \ii \mu(1 - \gamma) -\sum_{m=0}^n (m+1) \frac{\ii^{m}}{m!(2\mu)^{m}} \frac{(n+m)!}{(n-m)!} \Bigl(R_n\Bigl(\frac{i}{2 \mu}\Bigr)\Bigr)^{-1}.
\ee

 Setting $w = \frac{\ii}{ 2\mu} = \frac{ \ii \bar{\mu}}{2|\mu|^2},$ one  deduces     
\be \label{eq:7.9}
- C(n; \mu, \gamma) = \frac{1}{2 w}(1 - \gamma)   + (R_n(w))^{-1} \frac{d}{dw} (w R_n(w)).
\ee

Notice that $\Re \frac{\ii\bar{\mu}}{|\mu|^2} > 0$ implies $\re w > 0$, so we wish to prove that $C(n; \mu, \gamma) \neq 0$ for $\re w > 0.$
The case $n = 0$ is trivial because $\ii \mu (1- \gamma) - 1 \neq 0.$
We know that $ R_n(w)$ has no roots in the half plane $\re w \geq 0$. This implies that the roots of $w R_n(w) = 0$ lie in the half plane $\re w \leq 0.$ By the classical Gauss-Lucas theorem the roots of $\frac{d}{d w} (w R_n(w)) = 0$ lie in the convex set of the roots of $w R_n(w) = 0$ and one deduces that $\frac{d}{d w} (w R_n(w)) \neq 0$ for $\re w > 0.$ Thus for $\gamma = 1$ we have no eigenvalues of $G$.\\

For $0 < \gamma < 1$ we
 must examine the zeros of the function
 $$g_n(w) = \frac{1}{2w^2}(1- \gamma) + \frac{1}{w} + \sum_{j=1}^n \frac{1}{w - z_j},\: n \geq 1,$$
  where $z_j,\:\re z_j > 0,  j = 1,...,n,$ are the roots of $R_n(w) = 0$. 
We obtain
\be \label{eq:7bis}
\re g_n(w) = \frac {(1- \gamma)((\re w)^2 - (\im w)^2) + 2\re w |w|^2}{2 |w|^4} + \sum_{j=1}^n \frac{\re w - \re z_j}{|w- z_j|^2}.
\ee
 If $\re g_n(w) = 0$, we must have
$$2 \re w |w|^2 + (1- \gamma) ((\re w)^2 - (\im w)^2) < 0.$$
Setting $w = -\frac{1}{2 \lambda}$, the last inequality implies
$$|\re \lambda| < |\im \lambda|,\: \re \lambda > (1- \gamma)( (\re \lambda)^2- (\im \lambda)^2)$$
and we obtain that the eigenvalues of $G$ belong to the domain (\ref{eq:7.4}).

Passing to the case $\gamma > 1,$ we have the following
\begin{proposition}  For $\gamma > 1$ all eigenvalues $\lambda$ for which $(7.3)$ has a non trivial solution are real and they lie in the interval $[-\frac{1}{\gamma - 1}, -\infty).$
  Moreover, there is an infinite number  of  real eigenvalues of $G$.

\end{proposition}
{\it Proof.} 
To prove the existence of real eigenvalues, consider the polynomial
 $$F_n(w) = \Bigl[ \frac{1}{2}(1 - \gamma) + w\Bigr] R_n(w) + w^2 R_n'(w).$$
Clearly, $F_n(0) < 0$ and $F_n(w) \to +\infty$ as $w \to +\infty,$
so $F_n(w) = 0$ has at least one root $w_0$ in $\R^+$ and $C(n; \mu_0, \gamma) = 0$ for $\mu_0 = \frac{\ii}{2w_0}.$ \\

Now suppose that $w g_n(w_0) = 0,\: n \geq 1$ with $\re w_0 > 0, \: \im w_0 > 0.$ Then $\im g_n(w_0) = 0$ implies
\be \label{eq:7.10}
-\frac{(1- \gamma) \im w_0}{2 |w_0|^2} + \re w_0\Bigl[- \sum_{j=1}^n \frac{\im w_0}{|w_0- z_j|^2} + \sum_{j =1}^n \frac{ \im z_j}{|w_0- z_j|^2}\Bigr]\ee
$$ + \im w_0 \sum_{j=1}^n \frac{\re w_0 - \re z_j}{|w_0- z_j|^2} = 0.$$

 On the other hand, if $z_j$ with $\im z_j \neq 0$ is a root of $R_n(w)=0$, then $\bar{z_j}$ is also a root and
$$\frac{\im z_j}{|w_0 - z_j|^2} - \frac{\im z_j}{|w_0 - \bar{z}_j|^2}= \frac{\im z_j}{|w_0 - z_j|^2 |w_0 - \bar{z}_j|^2} \Bigl( |w_0- \bar{z}_j|^2 - |w_0 - z_j|^2\Bigr)$$
$$= \frac{4 \im w_0 (\im z_j)^2}{|w_0 - z_j|^2 |w_0 - \bar{z}_j|^2}.$$
Hence we can write (\ref{eq:7.10}) as follows
\be \label{eq:7.12}
\im w_0 \Bigl[ \frac{\gamma - 1}{2|w_0|^2} - \sum_{j=1}^n \frac{\re z_j}{|w_0- z_j|^2} + \sum_{\im z_j > 0}\frac{4 \re w_0 (\im z_j)^2}{|w_0- z_j|^2 |w_0- \bar{z}_j|^2}\Bigr] = 0.
\ee
 The term in the brackets $[...]$ is positive, and one concludes that $\im w_0 = 0.$ The same argument works for $\gamma = 1$ since $z_j \neq 0.$ Thus for $\gamma = 1$ we may have only real roots and since $w R_n'(w) \neq 0$ for $w > 0$ we conclude that there are no roots of $g_n(w) = 0.$

 From $\re g_n(w_0) = 0$ , one deduces for the real roots $w_0$ the equality 
 $$\frac{1 - \gamma}{2 w_0} + 1 = \sum_{j=1}^n \frac{\re z_j}{|w_0 - z_j|^2} > 0$$
  and this yields for the eigenvalues $\lambda$ of $G$ the inequality
$$\lambda \leq -\frac{1}{\gamma - 1}.$$

It remains to show that we have an infinite number of real eigenvalues. It is not excluded that for $n \neq m$ the polynomials $F_n(w)$ and $F_m(w)$ have the same real positive root. If we assume that for $\re w > 0$ the sequence of polynomials $\{F_n(w)\}_{n=0}^{\infty}$ has only a finite number of real roots $w_1,...,w_N$, $w_j \in \R^+$, then there exists an infinite number of polynomials $F_{n_j}(w)$ having  the same root which implies that we have an eigenvalue of $G$ with infinite multiplicity. This is a contradiction, and the number of real eigenvalues of $G$ is infinite. \\

\begin{remark} With small modifications Propositions $7.1$ and $7.2$ can be established for the ball $\{x\in \R^d: \:|x| \leq 1\}$ and $d \geq 5$ odd, by using the modified Hankel functions
$$\frac{H_{n + d/2 -1}^{(1)}(r)}{r^{d/2 -1}}$$
and the eigenfunctions $Y_{n, m}(\omega)$ of the Laplace-Beltrami operator $-\Delta_{S^{d-1}}$ with eigenvalues $n(n +d-2)$.

\end{remark}

{\bf Acknowledgments.} Thanks are due to Georgi Vodev  for many useful discussions and remarks concerning the paper \cite{V1}. I am also grateful to Johannes Sj\"ostrand for the fruitful discussions on the results in Chapters 9 and 10 of his paper \cite{Sj}. 
Finally, I am grateful to the referees for their careful reading of the manuscript and for the valuable comments
and suggestions.



\begin{equation*}
\end{equation*}

\end{document}